\documentclass[11pt,a4paper,leqno,twoside, headinclude]{amsart}

\title{Mapping degrees of self-maps of simply-connected manifolds}
\def\titl{Mapping degrees of self-maps of simply-connected manifolds}
\def\auth{Manuel Amann}
\date{September 4th, 2011}

\thanks{The author was supported by a Research Grant of the German Research Foundation.}
\author{\auth}

\usepackage[english]{babel}
\usepackage[colorlinks,pdftex, plainpages=false]{hyperref}
\hypersetup{pdftitle=\titl, pdfauthor=\auth, pdftoolbar=false,
plainpages=false, hyperindex=true, pdfdisplaydoctitle=true}

\hypersetup{colorlinks,%
citecolor=black,%
filecolor=black,%
linkcolor=black,%
urlcolor=black,%
pdftex}

\usepackage{color}

\usepackage{amsmath, amssymb, amscd, amsthm}
\usepackage{stmaryrd}
\usepackage{latexsym}
\usepackage[all]{xy}
\usepackage{pb-diagram}
\usepackage{rotating}
\usepackage{multicol}
\usepackage{lscape}
\usepackage{wasysym}
\usepackage{longtable}
\xyoption{all}

\newtheorem{theo}{Theorem}[section]
\newtheorem{main}{Theorem}
\newtheorem*{main*}{Theorem}

\newtheorem{prop}[theo]{Proposition}
\newtheorem{defi2}[theo]{Definition}

\newtheorem{rem2}[theo]{Remark}
\newenvironment{rem}{\begin{rem2}\normalfont}{\hfill$\boxbox$\end{rem2}}
\newtheorem{lemma}[theo]{Lemma}
\newtheorem{cor}[theo]{Corollary}

\newtheorem*{conj*}{Conjecture}
\newtheorem*{theo*}{Theorem}
\newtheorem*{ques*}{Question}
\newtheorem*{mi2}{Main Idea}

\newtheorem{ex2}[theo]{Example}
\newenvironment{ex}{\begin{ex2}\normalfont}{\hfill$\boxbox$\end{ex2}}
\newtheorem{exer2}[theo]{Exercise}

\newtheorem{alg2}[theo]{Algorithm}

\newcommand{\cc}{{\mathbb{C}}}                                     % complex numbers
\newcommand{\nn}{{\mathbb{N}}}                                     % natural numbers
\newcommand{\qq}{{\mathbb{Q}}}                                     % rational numbers
\newcommand{\pp}{{\mathbf{P}}}                                     % polynomials, projective space
\newcommand{\s}{{\mathbb{S}}}                                      % sphere
\newcommand{\zz}{{\mathbb{Z}}}                                     % integers
\newcommand{\SU}{{\mathbf{SU}}}                                    % special unitary group
\newcommand{\Sp}{{\mathbf{Sp}}}                                    % symplectic group
\newcommand{\B}{{\mathbf{B}}}                                      % Lie group type B
\newcommand{\dif} {{\operatorname{d}}}                             % differential operator d
\newcommand{\In} {{\,\subseteq\,}}                                 % subset
\newcommand{\Hom}{{\operatorname{Hom}}}                            % homomorphisms
\newcommand{\Ext}{{\operatorname{Ext}}}                            % Ext-group
\newcommand{\id}{{\operatorname{id}}}                              % identity
\newcommand{\APL}{{\operatorname{A_{PL}}}}                         % polynomial differential forms
\newcommand{\ADR}{{\operatorname{A_{DR}}}}                         % differential forms
\newcommand{\vproof}{{\begin{flushright} \qed                      % empty proof
                      \end{flushright}}}

\newcommand{\comment}[1]{}                                         % insert a large comment
\newcommand{\xto}[1]{\xrightarrow{#1}}                             % abbreviation for \xrightarrow
\newcommand{\hto}[1]{\overset{#1}{\hookrightarrow}}                % abbreviation for labelled \hookrightarrow
\newcommand{\biq}[2]{#1\;\!\!\!\sslash \;\!\!\!#2}                 % creates a biquotient
\newcommand{\ack}{\noindent\textbf{Acknowledgements. }}            % Acknowledgements
\newcommand{\str}{\noindent\textbf{Structure of the article. }}    % Structure of the article
\newenvironment{prf}{\begin{proof}[\textsc{Proof}]} {\end{proof}}     % alternative environment proof
\begin{document}

\maketitle \thispagestyle{empty}

%%%%%%%%%%%%%%%%%%%%%%%%%%%%%%%%%% Abstract%%%%%%%%%%%%%%%%%%%%%%%%%%%%%%%%%%%%%%%%

\begin{abstract}
An oriented compact closed manifold is called \emph{inflexible} if the set of mapping degrees ranging over all continuous self-maps is finite.
Inflexible manifolds have become of importance in the theory of functorial semi-norms on homology.
Although inflexibility should be a generic property in large dimensions, not many simply-connected examples are known. We show that from a certain dimension on there are infinitely many inflexible manifolds in each dimension. Besides, we prove \emph{flexibility} for large classes of manifolds and, in particular, as a spin-off, for homogeneous spaces. This is an outcome of a lifting result which also permits to generalise a conjecture of Copeland--Shar to the ``real world''.

Moreover, we then provide examples of simply-connected smooth compact closed manifolds in each dimension from dimension $70$ on which have the following properties: They do not admit any self-map which reverses orientation. (For this we consider the lack of orientation reversal in the strongest sense possible, i.e.~we prove the non-existence of any self-map of \emph{arbitrary negative} mapping degree.) Moreover, the manifolds neither split as non-trivial Cartesian products nor as non-trivial connected sums.
\end{abstract}

%%%%%%%%%%%%%%%%%%%%%%%%%%%%%%%%%% Introduction %%%%%%%%%%%%%%%%%%%%%%%%%%%%%%%%%%%

\section*{Introduction}
The \emph{mapping degree} of a continuous map is a fundamental invariant in manifold topology. For a continuous map $f: M\to N$ between oriented connected compact closed $n$-manifolds $M$ and $N$ with fundamental classes/orientation classes $[M]$ and $[N]$ it is defined by
\begin{align*}
\deg f&= d \intertext{with $d\in \zz$ being uniquely determined by}
H_n(f)([M])&=d\cdot[N]
\end{align*}
Here
\begin{align*}
H_n(f): H_n(M,\zz)\to H_n(N,\zz)
\end{align*}
denotes the induced map in singular homology. As it becomes obvious in the case of $M=N=\s^1$ for example, the mapping degree expresses in some sense how often the first manifold ``winds around'' its target by means of $f$.

In this article we are concerned with the question of how generic special mapping degrees are in large dimensions. (In small dimensions no such genericity can be expected, of course.)

Fix an oriented compact closed manifold $M$.
We shall focus on two aspects of being ``special'', which meet large interest. We shall consider manifolds $M$ for which either
\begin{itemize}
\item the set $\{\deg f\mid f:M\to M\}$ is finite or for which
\item there is no map $f: M\to M$ with $\deg f<0$.
\end{itemize}
By composition of self-maps it is obvious that the first case is equivalent to demanding that
\begin{align*}
\{\deg f\mid f:M\to M\}\In\{-1,0,+1\}
\end{align*}
In \cite{CL11} these manifolds were named \emph{inflexible} manifolds. In that article these manifolds were used to construct examples of certain functorial semi-norms on singular homology---cf.~remark \ref{rem03}. The examples of inflexible manifolds used there are adaptations of examples given in the article \cite{AL00}. Here the set of all homotopy classes $[M,N]$ of homomorphisms between minimal differential graded algebras $M$ and $N$ is studied extensively via a newly developed obstruction theory and, in particular, examples of minimal algebras which only admit trivial endomorphisms are constructed. In the recently appeared article \cite{CV11}---presumably and reportedly still under revision---theorem \cite{CV11}.1.5., p.~2, yields that every finite connected graph with trivial automorphism group leads to an elliptic space of even dimension without non-trivial self-homotopy equivalences. Examples of such graphs---like the following one---however, are well-known and easy to find:
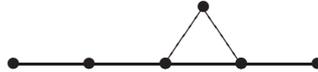
\begin{align*}
\xy<1cm, 0cm>:
(0,0)*=0{\bullet}="0";
(1,0)*=0{\bullet}="1" **@{-};
(2,0)*=0{\bullet}="2" **@{-};
(3,0)*=0{\bullet}="3" **@{-};
(4,0)*=0{\bullet}="4" **@{-};
(3,0)*=0{\bullet}="3" **@{-};
(2.5,0.75)*=0{\bullet}="5" **@{-};
(2,0)*=0{\bullet}="2" **@{-};
\endxy
\end{align*}

\vspace{5mm}

As for the second case we remark that a map $f$ with $\deg f<0$ clearly reverses the orientation of the manifold. The question of orientation reversal is rather prominent in topology. Different notions of this concept varying the categories in question and the definition of orientation reversal itself can be found. The question whether certain orientation-reversing maps on manifolds exist or not is rather classical and difficult---cf.~\cite{Pup95} for example. We shall exclusively work with continuous maps and whenever we speak of orientation reversal we shall mean orientation reversal in the weakest sense possible, i.e.~the existence of a self-map of \emph{arbitrary negative} mapping degree.

It seems to be the case that the question of inflexibility as well as the problem of lack of orientation reversal becomes much harder when focussing on the simply-connected case, since the fundamental group and its endomorphisms do pose certain obstructions---cf.~\cite{CL11}, \cite{Mue09}.
In \cite{Mue09} manifolds which do not possess a self-map of mapping degree $-1$ were called \emph{strongly chiral}. In this article strongly chiral manifolds in every dimension starting in dimension three, respectively in dimension seven for the simply-connected case were constructed. The constructions rely on torsion obstructions and extensively use Cartesian products in order to find higher dimensional examples.

Although the simply-connected case makes life much harder it does provide one advantage: Rational methods apply, i.e.~Rational Homotopy Theory enters the stage and does provide an efficient tool to study the mapping degree. We want to point out that one advantage of the rational approach is that it will enable us to consider arbitrary negative mapping degree, not only degree $-1$.

We consider the morphisms of commutative differential graded algebras induced by $f$ on smooth respectively polynomial differential forms:
\begin{align*}
\APL(f)&:~ \APL(N)\to \APL(M)\\
\ADR(f)&:~ \ADR(N)\to \ADR(M)\\
\end{align*}
Choosing representatives of the duals of fundamental classes $[M]$ and $[N]$ we may define the mapping degrees of $\APL(f)$ and $\ADR(f)$ via their induced maps in (co)homology $H(\APL(f))$ and $H(\ADR(f))$ as above. By a DeRham theorem this mapping degree coincides with the one defined above.

Thus we shall exclusively deal with simply-connected spaces and we shall draw on the rational approach. In analogy to the definition above it is hence possible to speak of inflexible commutative differential graded algebras or of orientation reversal whenever the algebras satisfy Poincar\'e Duality. We shall follow the textbook \cite{FHT01} in notation and we recommend it as a reference for the rational methods applied.

\vspace{5mm}

Coming back to our two main problems, we remark that in both cases several examples of manifolds satisfying the respective conditions are known---cf.~\cite{CL11}, \cite{Mue09}. However, still many open questions remain. For example, in the case of inflexible manifolds no examples of dimensions not divisible by $4$ were known; large gaps between the dimensions of the known examples exist. In both cases, examples in varying dimensions are constructed mainly using Cartesian products and partly connected sums. Moreover, as for strongly chiral manifolds torsion invariants play a crucial role. Our rational approach thus even permits to widen the notion of non-existence of orientation reversal as indicated.

We call a simply-connected manifold \emph{irreducible} if it does not split as a non-trivial Cartesian product of connected compact closed manifolds. We call it \emph{prime} if it does not split non-trivially as a connected sum. The following theorem closes existing gaps.
\begin{main}\label{theoA}
In each dimension $231+4i$ (for $i\geq 0$) there is an inflexible irreducible prime simply-connected smooth compact closed manifold $M_i$.

In each dimension greater than or equal to $921$ there are infinitely many simply-connected inflexible compact closed smooth manifolds.
\end{main}
\begin{prf}
The first assertion follows from a combination of lemma \ref{lemma01}, theorem \ref{theo01}, lemma \ref{lemma02} and remark \ref{rem02}. The second assertion is corollary \ref{cor01}.
\end{prf}

Indeed, inflexibility is expected to hold in the generic case. This is why it is certainly desirable to identify large classes of flexible spaces. Using obstruction theory we shall provide an approach to this which allows to relate the existence of a map on a rationalisation of a manifold in question to the existence of suitable self-maps on the actual manifold itself. As a corollary of a far more general theorem (which also generalises the flexibility in the formal case) we obtain
\begin{main}\label{theoC}
Simply-connected biquotients are flexible. Besides, so are oriented homogeneous spaces with connected denominator group.
\end{main}
\begin{prf}
This is the content of corollary \ref{cor01}.
\end{prf}
A conjecture of Copeland--Shar claims that there are either infinitely many homotopy classes of maps between certain rational spaces $M$, $N$ or the set $[M,N]$ is trivial. This conjecture was shown to hold true on large classes of spaces, but to be false in larger generality---cf.~\cite{AL00}. As another outcome of our approach, we provide a generalisation of this conjecture to not necessarily rational spaces $M$, $N$ of the rational homotopy of the same large class of spaces.

As for the lack of orientation reversal we are not only able to exclude the possibility of self-maps of arbitrary negative degree, but the examples we provide also neither split non-trivially as products nor as connected sums.
\begin{main}\label{theoB}
In every dimension $\dim M\geq 70$ there is an irreducible prime simply-connected smooth compact closed manifold which does not possess any self-map of negative mapping degree.
\end{main}

\vspace{5mm}

\textbf{As a general convention cohomology will be taken with rational coefficients unless stated differently.}

\vspace{5mm}

\str In section \ref{sec01} we provide a sequence of infinitely many inflexible manifolds. Indeed, these manifolds will cover all dimensions congruent to $3$ modulo $4$ from dimension $231$ on. They constitute the first examples of such manifolds in odd dimensions. We shall make use of these examples to prove the main ingredients to theorem \ref{theoA}.

In section \ref{sec02} we provide a theorem which grants flexibility of a manifold under mild assumptions on a map on its rationalisation. This will allow us to prove that whenever the differential of the minimal model of a space is well-behaved with respect to a lower grading, the space will be flexible. This proves flexibility for simply-connected two-stage spaces and theorem \ref{theoC}, in particular. Moreover, we deal with the conjecture by Copeland--Shar.

In section \ref{sec04} irreducible examples of manifolds without orientation reversal are constructed. This leads to the proof of theorem \ref{theoB}.

\vspace{5mm}

\ack The author is very grateful to Clara L\"oh for several helpful remarks.

%%%%%%%%%%%%%%%%%%%%%%%%%%%%%%%%%% Section 1 %%%%%%%%%%%%%%%%%%%%%%%%%%%%%%%%%%%%%%

\section{Irreducible inflexible manifolds in infinitely many odd dimensions}\label{sec01}
We define the sequence of commutative differential graded algebras
\begin{align*}
\big((A_i,\dif) \mid i\geq 0\big)
\end{align*}
by $A_i=\bigwedge V_i$ where $V_i$ is the graded rational vector space generated by the elements
\begin{align*}
&x_1,x_2,y_1,y_2,y_3,z,z'
\intertext{of degrees}
&\deg x_1=4, \deg x_2=6, \deg y_1=27, \deg y_2=29, \deg y_3=31,\\& \deg z=77, \deg z'=75+4i
\end{align*}
The differential $\dif$ is defined by
\begin{align*}
&\dif x_1=\dif x_2=0, \dif y_1=x_1^4x_2^2, \dif y_2=x_1^3x_2^3, \dif y_3=x_1^2x_2^4, \\&\dif z=x_1x_2^3y_1y_2-x_1^2x_2^2y_1y_3+x_1^3x_2y_2y_3+x_2x_1^{18}+x_2^{13}, \\&\dif z'=x_1^{19+i}
\end{align*}
\begin{lemma}\label{lemma01}
Suppose that $i\geq 0$.
\begin{itemize}
\item The algebra $(A_i,\dif)$ is an elliptic minimal Sullivan algebra and
\item its formal dimension satisfies $\dim (A_i,\dif)=231+4i$.
\item A volume form is given by the class of $x_2^{26}z'-x_1^{15+i}x_2^{24}y_1$.
\item The algebra $(A_i,\dif)$ can be realised as the minimal model of an irreducible compact closed simply-connected smooth manifold $M_i$ of the same odd dimension $\dim M_i=231+4i$.
\end{itemize}
\end{lemma}
\begin{prf}
It is obvious that $\dif^2=0$ on $A_i$, since
\begin{align*}
\dif^2 z&=\dif \big(x_1x_2^3y_1y_2-x_1^2x_2^2y_1y_3+x_1^3x_2y_2y_3+x_2x_1^{18}+x_2^{13}\big)
\\&= x_1^5x_2^5y_2-x_1^4x_2^6y_1-x_1^6x_2^4y_3+x_1^4x_2^6y_1+x_1^6x_2^4y_3-x_1^5x_2^5y_2
\\&=0
\end{align*}
Also minimality is clear. Let us now see that the algebra is elliptic. As it is finitely generated it suffices to see that it has finite dimensional homology. For this we prove that every cohomology class is nilpotent. It suffices to see that a certain power of $[x_1]$ and of $[x_2]$ vanishes. Due to $\dif z'=x_1^{19+i}$ this follows from
\begin{align*}
[x_2]^{26}&=[\dif z- x_1x_2^3y_1y_2+x_1^2x_2^2y_1y_3-x_1^3x_2y_2y_3-x_2x_1^{18}]^2\\
&=-2[- x_1^{19}x_2^4y_1y_2+x_1^{20}x_2^3y_1y_3-x_1^{21}x_2^2y_2y_3]+[x_2^2x_1^{36}]
\\&=2[x_1^{17}\dif(y_1y_2y_3)]+[\dif (x_1^{32}y_1)]
\\&=0
\end{align*}
(For this note that $(- x_1x_2^3y_1y_2+x_1^2x_2^2y_1y_3-x_1^3x_2y_2y_3)^2$ vanishes.)

\vspace{5mm}

As the algebra is elliptic we may compute its formal dimension by the formula (cf.~\cite{FHT01}.32, p.~434)
\begin{align*}
\dim (A_i,\dif)=&(\deg y_1+\deg y_2 +\deg y_3 + \deg z + \deg z')\\& - (\deg x_1-1+\deg x_2-1)\\=&231+4i
\end{align*}

\vspace{5mm}

Let us now see that $x_2^{26}z'-x_1^{15+i}x_2^{24}y_1$ (of degree $231+4i$) represents a non-vanishing cohomology class.
First of all we see that
\begin{align*}
\dif\big(x_2^{26}z'-x_1^{15+i}x_2^{24}y_1\big)=x_1^{19+i}x_2^{26}-x_1^{15+i}x_2^{24}x_1^4x_2^2=0
\end{align*}

The form is not exact, however: Assume the contrary, i.e.~the existence of a class $a\in (A^{230+41}_i,\dif)$ with $\dif a=x_2^{26}z'-x_1^{15+i}x_2^{24}y_1$. Then $a$ is necessarily of the form
\begin{align*}
&p_1y_1y_2+p_2y_1y_3+p_3y_2y_3+p_4y_1z+p_5y_2z+p_6y_4z+p_7y_1z'+p_8y_2z'
\\&+p_9y_3z'+p_{10}zz'+p_{11}y_1y_2y_3z+p_{12}y_1y_2y_3z'+p_{13}y_1y_2zz'
\\&+p_{14}y_1y_3zz'+p_{15}y_2y_3zz'
\end{align*}
with homogeneous polynomials $p_i\in \qq[x_1,x_2]$. We compute $\dif a$ as
\begin{align*}
&p_1x_1^4x_2^2y_2-p_1y_1x_1^2x_2^2+p_2x_1^4x_2^2y_3-p_2y_1x_1^2x_2^4
\\+&p_3x_1^3x_2^3y_3-p_3y_2x_1^2x_2^4+p_4x_1^4x_2^2z+p_5x_1^3x_2^3z+p_6x_1^2x_2^4z
\\-&(p_4y_1+p_5y_2+p_6y_3)(x_1x_2^3y_1y_2-x_1^2x_2^2y_1y_3+x_1^3x_2y_2y_3+x_2x_1^{18}+x_2^{13})
\\+&p_7x_1^4x_2^2z'-p_7y_1x_1^{19+i}+p_8x_1^3x_2^3z'-p_8y_2x_1^{19+i}+p_9x_1^2x_2^4z'-p_9y_3x_1^{19+i}
\\+&p_{10}(x_1x_2^3y_1y_2-x_1^2x_2^2y_1y_3+x_1^3x_2y_2y_3+x_2x_1^{18}+x_2^{13})z'-p_{10}zx_1^{19+i}
\\+&p_{11}x_1^4x_2^2y_2y_3z-p_{11}y_1x_1^3x_2^3y_3z+p_{11}y_1y_2x_1^2x_2^4z-p_{11}y_1y_2y_3(x_2x_1^{18}+x_2^{13})
\\+&p_{12}x_1^4x_2^2y_2y_3z'-p_{12}y_1x_1^3x_2^3y_3z'+p_{12}y_1y_2x_1^2x_2^4z'-p_{12}y_1y_2y_3x_1^{19+i}
\\+&p_{13}x_1^4x_2^2y_2zz'-p_{13}y_1x_1^3x_2^3zz'+p_{13}y_1y_2(x_2x_1^{18}+x_2^{13})z'-p_{13}y_1y_2zx_1^{19+i}
\\+&p_{14}x_1^4x_2^2y_3zz'-p_{14}y_1x_1^2x_2^4zz'+p_{14}y_1y_3(x_2x_1^{18}+x_2^{13})z'-p_{14}y_1y_3zx_1^{19+i}
\\+&p_{15}x_1^3x_2^3y_3zz'-p_{15}y_2x_1^2x_2^4zz'+p_{15}y_2y_3(x_2x_1^{18}+x_2^{13})z'-p_{15}y_2y_3zx_1^{19+i}
\end{align*}
Comparing coefficients in $z'$ we obtain the equation
\begin{align*}
&x_2^{26}=\\
&p_7x_1^4x_2^2+p_8x_1^3x_2^3+p_9x_1^2x_2^4
\\&+p_{10}(x_1x_2^3y_1y_2-x_1^2x_2^2y_1y_3+x_1^3x_2y_2y_3+x_2x_1^{18}+x_2^{13})
\\&+p_{12}x_1^3x_2^2y_2y_3-p_{12}y_1x_1^3x_2^3y_3+p_{12}y_1y_2x_1^2x_2^4+p_{13}x_1^4x_2^2y_2z-p_{13}y_1x_1^3x_2^3z
\\&+p_{13}y_1y_2(x_2x_1^{18}+x_2^{13})+p_{14}x_1^4x_2^2y_3z-p_{14}y_1x_1^2x_2^4z+p_{14}y_1y_3(x_2x_1^{18}+x_2^{13})
\\&+p_{15}x_1^3x_2^3y_3z-p_{15}y_2x_1^2x_2^4z+p_{15}y_2y_3(x_2x_1^{18}+x_2^{13})
\end{align*}
This lets us conclude that $p_{10}=x_2^{13}$. Comparing coefficients in $y_2y_3z'$ we obtain the equation
\begin{align*}
p_{12}x_1^4x_2^2+p_{15}(x_2x_1^{18}+x_2^{13})=-p_{10}x_1^3x_2
\end{align*}
Since $p_{10}x_1^3x_2=x_1^3x_2^{14}$ we deduce that $p_{15}=-x_2x_1^3$.

Stepping back to the expression for $\dif a$ we use that there is no summand with a factor $z$ in our designated volume form. More precisely, we compare coefficients in $y_2y_3z$, which yields the equation
\begin{align*}
p_{11}x_1^4x_2^2=p_{15}x_1^{19+i} \Leftrightarrow p_{11}x_1^4x_2^2=-x_2x_1^{22+i}
\end{align*}
This is a clear contradiction, as such a $p_{11}\in \qq[x_1,x_2]$ cannot exist.

\vspace{5mm}

Finally, we apply theorem \cite{Sul77}.13.2, p.~321, which allows us to realise the algebras $(A_i,\dif)$ in the desired way:
Since the $(A_i,\dif)$ are simply-connected and elliptic, their cohomology algebras satisfy Poincar\'e duality---cf.~\cite{FHT91}, theorem A, p.~70. Thus it suffices to observe that $\dim (A_i,\dif)\not \equiv 0 \mod 4$ in order to realise the algebras. The manifolds are irreducible, since the minimal models do not split non-trivially as products.
\end{prf}

\begin{theo}\label{theo01}
For $i\geq 0$ the manifolds $M_i$ are inflexible.
\end{theo}
\begin{prf}
We prove that the algebras $(A_i,\dif)$ are inflexible, which requires the manifolds $M_i$ to be so either. This is due to the fact that every continuous map between simply-connected topological spaces has a Sullivan representative---cf.~\cite{FHT01}.12.9, p.~153.

Suppose there is a morphism of differential graded algebras
\begin{align*}
f: (A_i,\dif)\to (A_i,\dif)
\end{align*}
By degree we obtain that
\begin{align*}
f(x_1)=&k_1 x_1\\
f(x_2)=&k_2 x_2\\
f(y_1)=&k_3 y_1\\
f(y_2)=&k_4 y_2\\
f(y_3)=&k_5 y_3+k_6x_1y_1\\
f(z)=&k_7 z + p_1 y_1+p_2y_2+p_3y_3\\
f(z')=&k_8 z' + p_4y_1+p_5y_2+p_6y_3+p_7y_1y_2y_3+p_8z+p_9y_1y_2z+p_{10}y_1y_3z\\&+p_{11}y_2y_3z
\end{align*}
for $k_i\in \qq$ and $p_i\in \qq[x_1,x_2]$. Since $\dif (x_1y_1)=x_1^5x_2^2\not \in \langle x_1^2x_2^4\rangle_\qq$, we obtain $k_6=0$, since $f$ commutes with differentials. Due to the same reasoning the terms
\begin{align*}
f(\dif y_1)&=f(x_1^4x_2^2)=k_1^4k_2^2x_1^4x_2^2
\\\dif(f(y_1))&=k_3 \dif y_1=k_3 x_1^4x_2^2
\end{align*}
coincide. Thus we obtain that $k_3=k_1^4k_2^2$. In the same vein we compute $k_4=k_1^3k_2^3$ and $k_5=k_1^2k_2^4$.

Since $f\circ \dif=\dif \circ f$ we compute that
\begin{align*}
f(\dif z)=&f\big(x_1x_2^3y_1y_2-x_1^2x_2^2y_1y_3+x_1^3x_2y_2y_3+x_2x_1^{18}+x_2^{13}\big)\\=&k_1^8k_2^8 x_1x_2^3y_1y_2-k_1^8k_2^8x_1^2x_2^2y_1y_3+k_1^8k_2^8x_2y_2y_3+k_1^{18}k_2x_2x_1^{18}\\&+k_2^{13}x_2^{13}
\intertext{equals}
\dif(f(z))=&k_7(x_1x_2^3y_1y_2-x_1^2x_2^2y_1y_3+x_1^3x_2y_2y_3+x_2x_1^{18}+x_2^{13})\\&+p_1x_1^4x_2^2+p_2x_1^3x_2^3+p_3x_1^2x_2^4
\end{align*}
A comparison of coefficients yields
\begin{align*}
k_7&=k_1^8k_2^8=k_1^{18}k_2=k_2^{13}
\intertext{and}
\dif(p_1y_1+p_2y_2+p_3y_3)&=p_1x_1^4x_2^2+p_2x_1^3x_2^3+p_3x_1^2x_2^4=0
\end{align*}
Consequently, we have that $k_1^{\frac{20}{21}}=1$ unless $k_1=0$ or $k_2=0$ and we derive that either $k_2=0$ or $k_1=\pm 1 \wedge k_2=1$.

Let us discuss these two cases separately. Suppose first that $k_2=0$. We compute
\begin{align*}
f\big(x_2^{26}z'-x_1^{15+i}x_2^{24}y_1\big)=k_2^{26}x_2^{26}f(z')-k_1^{19+i}k_2^{26}x_1^{15+i}x_2^{24}y_1=0
\end{align*}
By lemma \ref{lemma01}, the form $x_2^{26}z'-x_1^{15+i}x_2^{24}y_1$ is a volume form. Thus the mapping degree of $f$ is $0$, i.e.~the manifold $M_i$ is inflexible.

We shall do similar in the case $k_1=\pm 1 \wedge k_2=1$. For this we need to investigate the image of $z'$ under $f$. We compute that
\begin{align*}
\dif\big(f(z')\big)=&\dif(k_8 z' + p_4y_1+p_5y_2+p_6y_3+p_7y_1y_2y_3+p_8z+p_9y_1y_2z\\&+p_{10}y_1y_3z+p_{11}y_2y_3z)
\\=&k_8x_1^{19+i}+\dif T
\intertext{equals}
f(\dif z')=&(\pm 1)^{19+i}x_1^{19+i}
\end{align*}
where
\begin{align*}
T:=p_4y_1+p_5y_2+p_6y_3+p_7y_1y_2y_3+p_8z+p_9y_1y_2z+p_{10}y_1y_3z+p_{11}y_2y_3z
\end{align*}
and $\dif T$ is a sum of monomials which each have either $x_2$, $y_i$ or $z$ as a factor. A comparison of coefficients thus shows that $k_8=(\pm 1)^{19+i}$ and that $\dif T$ vanishes. Thus we have that $f(z')=(\pm 1)^{19+i}z'+T$ with $\dif T=0$. Thus $T$ defines a cohomology class.

Let us now use this information on the volume form $[x_2^{26}z'-x_1^{15+i}x_2^{24}y_1]$. We compute
\begin{align*}
&H(f)\big([x_2^{26}z'-x_1^{15+i}x_2^{24}y_1]\big)
\\=&[x_2^{26}((\pm 1)^{19+i}z'+T)-(\pm 1)^{15+i}x_1^{15+i}x_2^{24}y_1]
\\=&[(\pm 1)^{19+i}x_2^{26}z'-(\pm 1)^{15+i}x_1^{15+i}x_2^{24}y_1]+[x_2^{26}][T]
\\=&[(\pm 1)^{19+i}x_2^{26}z'-(\pm 1)^{15+i}x_1^{15+i}x_2^{24}y_1]
\\=&(\pm 1)^{i+1}[x_2^{26}z'-x_1^{15+i}x_2^{24}y_1]
\end{align*}
for the induced map in cohomology, since $[x_2^{26}]=0$---cf.~the proof of lemma \ref{lemma01}. Thus the mapping degree of $f$ is $(\pm 1)^{i+1}$ and $M_i$ is inflexible.
\end{prf}

\begin{theo}\label{theo06}
For $d_l\in \nn_0^k$ the Cartesian product
\begin{align*}
\prod_{1\leq l\leq k} M_{d_l}
\end{align*}
is inflexible.
\end{theo}
\begin{prf}
We prove the corresponding statement on minimal models, i.e.~the inflexibility of
\begin{align*}
\bigg(\bigotimes_{1\leq l\leq k} A_{d_l},\dif\bigg)
\end{align*}
where $\dif$ denotes the product differential.
Thus a volume form of the product is the product of the volume forms from lemma \ref{lemma01}.

Consider an endomorphism $f$ of this algebra. Let us first see that $x_{1,l}$---the element $x_1$ in component $l$---and $x_{2,l}$ have to map to a multiple of $x_{1,l'}$ respectively $x_{2,l'}$. By degree, we have
\begin{align*}
f(x_{1,l})&=\sum_j k_{1,j} x_{1,j}\\
f(x_{2,l})&=\sum_j k_{2,j} x_{2,j} \\
f(y_1,l)&=\sum_j k_{3,j} y_{1,j}
\end{align*}
Since $f$ commutes with the differential, we obtain that
\begin{align}
\label{eqn01}\dif(f(y_{1,l}))&=\dif \big(\sum k_{3,j} y_{1,j}\big)=\sum k_{3,j}x_{1,j}^4x_{2,j}^2 \intertext{equals}
\label{eqn02}f(\dif(y_{1,l}))&=f(x_{1,l}^4x_{2,l}^2)=\big(\sum k_{1,j} x_{1,j}\big)^4\big(\sum k_{2,j}x_{2,j}\big)^2
\\&\nonumber= \sum k_{1,j_1}^4k_{2,j_2}^2 x_{1,j_1}^4 x_{2,j_2}^2+\ldots
\end{align}
which implies that $k_{1,j_1}^4k_{2,j_2}=0$, i.e.~$k_{1,j_1}=0 \vee k_{2,j_2}=0$, for $j_1\neq j_2$. That is, we have
\begin{align}
\label{eqn03}k_{1,j_1}\neq 0 \Rightarrow k_{2,j_2}=0 \qquad \textrm{for $j_1\neq j_2$}
\intertext{and, conversely,}
\label{eqn04}k_{2,j_2}\neq 0 \Rightarrow k_{1,j_1}=0 \qquad \textrm{for $j_1\neq j_2$}
\end{align}
The morphism $f$ cannot map $x_{2,l}$ to zero unless its mapping degree is zero. Thus there is a $k_{2,j'}\neq 0$.

We shall also see that $f(x_{1,l})=0$ implies that $\deg f=0$. Indeed, in this case we obtain $f(\dif z'_l)=0$ and thus $f(z'_l)$ has to be a closed form; i.e.~a cohomology class. We shall now show that the cohomology class of the volume form of the product vanishes under $H(f)$, since the class of the factor belonging to the $l$-th component vanishes. For this we compute
\begin{align*}
f\big(x_{2,l}^{26}z_l'-x_{1,l}^{15+i}x_{2,l}^{24}y_{1,l}\big)=f\big(x_{2,l}^{26}\big)f\big(z_l'\big)
\end{align*}
as $f(x_{1,l})=0$. Since $f(z_l')$ is closed and since $x_{2,l}^{26}$ is exact---consequently, so is $f(x_{2,l}^{26})$---we obtain that $f(x_{2,l}^{26})f(z_l')$ is exact.

Thus we may assume that there is a $k_{1,j''}\neq 0$. With the implications \eqref{eqn03} and \eqref{eqn04} from above we directly see that $j'=j''$ and that all the other $k_{1,j}$ respectively all the other $k_{2,j}$ vanish. That is
\begin{align*}
f(x_{1,l})&=k_{1,j'} x_{1,j'}\intertext{and}
f(x_{2,l})&=k_{2,j'}x_{2,j'}
\intertext{In view of the relations \eqref{eqn01} and \eqref{eqn02} this implies that}
f(y_{1,l})&=k_{3,j'}y_{1,j'}
\end{align*}
with $k_{3,j'}=k_{1,j'}^4k_{2,j'}^2$.

In a similar way we derive that $f(z'_l)=k_{1,l}^{19+i}(z_{j'}'+T_{j'})$ (for $k_{1,l}^{19+i}\neq 0$ unless $\deg f=0$) with a closed form $T_{j'}\in \bigotimes_{1\leq l\leq k} A_{d_l}$. Thus for the volume form of the $l$-th factor we compute
\begin{align*}
f(x_{2,l}^{26}z_l'-x_{1,l}^{15+i}x_{2,l}^{24}y_{1,l})=k_{2,j'}^{26}k_{1,j'}^{19+i}x_{2,j'}^{26}(z_{j'}'+T_{j'})-k_{1,j'}^{19+i}k_{2,j'}^{26}x_{1,j'}^{15+i}x_{2,j'}^{24}y_{1,j'}
\end{align*}
Thus on cohomology we may compute the image of the class of the volume form of $\big(\bigotimes_{1\leq l\leq k} A_{d_l},\dif\big)$ as
\begin{align*}
&H(f)\big(\prod_l [x_{2,l}^{26}z_l'-x_{1,l}^{15+i}x_{2,l}^{24}y_{1,l}]\big)
\\=&
\prod_l \big(k_{2,\pi(l)}^{26}k_{1,\pi(l)}^{19+i}[x_{2,\pi(l)}^{26}z_{\pi(l)}'-x_{1,\pi(l)}^{15+i}x_{2,\pi(l)}^{24}y_{1,\pi(l)}]
\\&+k_{2,\pi(l)}^{26}k_{1,\pi(l)}^{19+i}[x_{2,\pi(l)}]^{26}[T]\big)
\\=&
\prod_l k_{2,\pi(l)}^{26}k_{1,\pi(l)}^{19+i}[x_{2,\pi(l)}^{26}z_{\pi(l)}'-x_{1,\pi(l)}^{15+i}x_{2,\pi(l)}^{24}y_{1,\pi(l)}]
\end{align*}
since $[x_{2,\pi(l)}]^{26}=0$. Here $\pi: \{1,\dots, k\}\to\{1,\dots, k\}$ denotes an arbitrary map. However, unless $\deg f=0$ it is now a direct observation that $\pi$ has to be a permutation. Denote by $|\pi|$ the order of $\pi$.

We project the morphism $f$ to a morphism on the $l$-th factor:
\begin{align*}
f_l&:   A_{d_l}\hto{} \bigg(\bigotimes_{1\leq m\leq k} A_{d_m}\bigg)\xto{f} \bigg(\bigotimes_{1\leq m\leq k} A_{d_m}\bigg) \to A_{d_l}
\end{align*}
Thus $f_l$ is an endomorphism of an inflexible algebra $A_{d_l}$---cf.~theorem \ref{theo01}.
Hence the same is true for $f_l^{|\pi|}$. Moreover, we have proved that $H(f^{|\pi|})=(H(f))^{|\pi|}$ maps the cohomology class of the chosen volume form of $A_{d_l}$ to a multiple of itself. Thus
\begin{align*}
(H(f))^{|\pi|}\big([x_{2,l}^{26}z_l'-x_{1,l}^{15+i}x_{2,l}^{24}y_{1,l}]\big)
&=(H(f_l))^{|\pi|}\big([x_{2,l}^{26}z_l'-x_{1,l}^{15+i}x_{2,l}^{24}y_{1,l}]\big)
\\&=\pm [x_{2,l}^{26}z_l'-x_{1,l}^{15+i}x_{2,l}^{24}y_{1,l}]
\end{align*}
unless the mapping degrees of both $f$ and $f_l$ vanish.

Thus we compute
\begin{align*}
&\pm\prod_l [x_{2,l}^{26}z_l'-x_{1,l}^{15+i}x_{2,l}^{24}y_{1,l}]
\\=&(H(f))^{|\pi|}\big(\prod_l [x_{2,l}^{26}z_l'-x_{1,l}^{15+i}x_{2,l}^{24}y_{1,l}]\big)
\\=&\prod_l \big(k_{2,\pi(l)}^{26}k_{1,\pi(l)}^{19+i}\big)^{|\pi|} [x_{2,l}^{26}z_l'-x_{1,l}^{15+i}x_{2,l}^{24}y_{1,l}]
\end{align*}
Consequently, we obtain that
\begin{align*}
\prod_l \big(k_{2,\pi(l)}^{26}k_{1,\pi(l)}^{19+i}\big)^{|\pi|}&=\pm 1\intertext{respectively}
\bigg|\prod_l k_{2,\pi(l)}^{26}k_{1,\pi(l)}^{19+i}\bigg|&=1
\intertext{i.e.}
H(f)\big(\prod_l [x_{2,l}^{26}z_l'-x_{1,l}^{15+i}x_{2,l}^{24}y_{1,l}]\big)&=\pm \prod_l [x_{2,l}^{26}z_l'-x_{1,l}^{15+i}x_{2,l}^{24}y_{1,l}]
\end{align*}
unless $\deg f=0$. This proves that the product $\big(\bigotimes_{1\leq l\leq k} A_{d_l},\dif\big)$ is inflexible.
\end{prf}

\begin{cor}
In each dimension greater than or equal to $921$ there are infinitely many simply-connected inflexible smooth compact closed manifolds.
\end{cor}
\begin{prf}
Consider the manifolds
\begin{align*}
M_i, M_0\times M_i, M_0^2\times M_i, M_0^3\times M_i
\end{align*}
to get inflexible examples of dimensions
\begin{align*}
231+4i, 462+4i, 693+4i, 924+4i
\end{align*}
according to theorem \ref{theo06} and lemma \ref{lemma01}. This covers all congruence classes modulo $4$ from dimension $921$ on.

In order to get infinitely many examples we apply theorem \cite{CL11}.9.1 to these manifolds. We obtain that
\begin{align*}
M_i^{\#k}, (M_0\times M_i)^{\#k}, (M_0^2\times M_i)^{\#k}, (M_0^3\times M_i)^{\#k}
\end{align*}
are inflexible. (For this it is easily checked that the rational homotopy group below top dimension vanishes in each respective case; indeed, this holds true, since it holds true on each $M_i$ already.)
\end{prf}

\begin{rem}
From the proof of theorem \ref{theo01} we observe that the manifolds we constructed there do not even possess an orientation reversing self-map if $i$ is odd, i.e.~in dimensions $235+8i'$ for $i'\in \nn_0$. The same is true for direct products of the manifolds $M_{i'}$ of pairwise different dimensions, as can be derived from the proof of theorem \ref{theo06}. (Indeed, in this case the permutation $\pi$ has to be the identity, i.e~$|\pi|=1$, as we obtain that $f(z_j')\in \langle z_j'\rangle \oplus \ker \dif$ for degree reasons.) We refer to section \ref{sec04} for further examples of this.
\end{rem}

\begin{rem}\label{rem03}
As has been done in theorem \cite{CL11}.1.2, p.~2, one may use the existence of inflexible manifolds $M$ to derive properties of functorial semi-norms on singular homology. In the theorem the existence of functorial semi-norms on singular homology which assume only positive, finite values on certain simply-connected manifolds of dimension divisible by $4$ is proved---cf.~corollary \cite{CL11}.7.3.

For this the \emph{domination semi-norm} (cf.~definition \cite{CL11}.7.1, p.~21) associated with a manifold is defined. By proposition \cite{CL11}.7.2, p.~21, and the proof of corollary \cite{CL11}.7.3, p.21, there is a positive, finite functorial semi-norm on singular homology of degree $\dim M$ associated to the domination semi-norm---cf.~theorem \cite{CL11}.4.2, p.~8---which takes the value $1$ on the fundamental class of $M$.

The fact that it only assumes the value $1$ is now easily seen to be indebted to $M$ being inflexible. Indeed, the methods in \cite{CL11} generalise to our examples of inflexible manifolds thereby proving the existence of functorial semi-norms on singular cohomology which assume only finite, positive values on certain simply-connected manifolds the dimension of which may also be odd or congruent to $2$ modulo $4$.
\end{rem}

We end this section by providing a method and another concrete example of how to construct further irreducible manifolds in even dimensions not necessarily divisible by four.
\begin{ex}\label{ex01}
Using the inflexible examples from \cite{CL11}.8, p.~23, it is rather easy to construct further even-dimensional simply-connected inflexible manifolds. Indeed, one may even produce examples the dimensions of which are congruent to two modulo four. Let us give one explicit example of this form---indeed, this will be the most complicated example.

We consider the algebra from example \cite{CL11}.8.4, p.~24, which is given by $V=\langle x_1,x_2,y_1,y_2,y_3,z\rangle$ with $\deg x_1=2, \deg x_2=4, \deg y_1=9, \deg y_2=11, \deg y_3=13, \deg z=35$ and with
\begin{align*}
\dif x_1=\dif x_2&=0
\\\dif y_1&=x_1^3x_2
\\\dif y_2&=x_1^2x_2^2
\\\dif y_3&=x_1x_2^3
\\\dif z&=x_2^4y_1y_2-x_1x_2^3y_1y_3+x_1^2x_2^2y_2y_3+x_1^{18}+x_2^{9}
\end{align*}
Let us form a rational fibration with $(\bigwedge V,\dif)$ as a basis and fibre (a rational) $\s^2$:
\begin{align*}
\big(\bigwedge V \otimes \bigwedge \langle \bar x_2, x_2'\rangle,\dif\big)
\end{align*}
defined by $\dif \bar x_2=0$, $\dif x_2'=\bar x_2^2-x_2$. It is then obvious that a minimal model for the total space is given by
$V'=\langle x_1, \bar x_2,y_1,y_2,y_3,z\rangle$ with $\deg x_1=2, \deg \bar x_2=2, \deg y_1=9, \deg y_2=11, \deg y_3=13, \deg z=35$ and with
\begin{align*}
\dif x_1=\dif \bar x_2&=0
\\\dif y_1&=x_1^{3}\bar x_2^2
\\\dif y_2&=x_1^{2}\bar x_2^{4}
\\\dif y_3&=x_1\bar x_2^{6}
\\\dif z&=\bar x_2^{8}y_1y_2-x_1\bar x_2^{6}y_1y_3+x_1^2\bar x_2^4y_2y_3+x_1^{18}+\bar x_2^{18}
\end{align*}
The Sullivan algebra $(\bigwedge V',\dif)$ consequently is elliptic of dimension $66 \equiv 2 \mod 4$. As $x_2^{16}$ represents a volume form of $(\bigwedge V,\dif)$---cf.~proposition \cite{CL11}.8.6, p.~25, a volume form of $(\bigwedge V',\dif)$ is represented by $\bar x_2^{33}$. As $\dim (\bigwedge V',\dif)\equiv 2 \mod 4$, this algebra can be realised by a simply-connected smooth compact closed manifold $M^{66}$. Let us now see that it is inflexible.

In analogy to the arguments given in theorem \ref{theo01} we proceed as follows:
It is easy to see that an arbitrary endomorphism $f: (\bigwedge V',\dif)\to \bigwedge V',\dif)$ has to satisfy
\begin{align*}
f(x_1)&=k_1 x_1\\
f(\bar x_2)&=k_2 \bar x_2\\
f(y_1)&=k_1^{3}k_2^2 y_1\\
f(y_2)&=k_1^{2}k_2^4 y_2\\
f(y_3)&=k_1k_2^6 y_3+T_1\\
f(z)&=k_3 z+T_2
\end{align*}
with rational $k_i$ and closed forms $T_1, T_2\in \bigwedge V'$. The differential of $z$ then yields the relations
\begin{align*}
k_1^{5}k_2^{14}=k_1^{18}=k_2^{18}
\end{align*}
This amounts to the contradiction $k_2^{\frac{14}{13}}=1$ unless $k_2\in\{0,1,-1\}$. Thus we compute the image of the volume form as
\begin{align*}
f([\bar x_2^{33}])=k_2^{33}[x_2]^{33}\in\{\pm [x_2]^{33}, 0\}
\end{align*}
Consequently, the manifold $M$ is inflexible.

One may proceed similarly on the other examples from \cite{CL11}.8, p.~23, again using rational fibrations with fibres given by $\big(\bigwedge\langle \bar x, x' \rangle,\bar \dif\big)$ with $\deg x$ even, $\bar \dif \bar x=0$, $\bar \dif x'=\bar x^l$. The differential $\dif$ in the model of the fibration then is such that $\bar \dif x'=\bar x^l-x$ (where $x\in\{x_1,x_2\}$ is chosen according to whether a power of $x_1$ or of $x_2$ represents a volume form). This produces a lot of examples in even dimensions.
\end{ex}

%%%%%%%%%%%%%%%%%%%%%%%%%%%%%%%%%% Section 2 %%%%%%%%%%%%%%%%%%%%%%%%%%%%%%%%%%%%%%

\section{From rational to real}\label{sec02}

As it is likely that in the generic case a manifold is inflexible, it is certainly desirable to give sufficient criteria for flexibility.
For this we shall first prove a result which helps us relate the flexibility of a space to the one of its rationalisation---thus generalising results in \cite{Shi79}. Then we prove flexibility for several large classes of manifolds amongst which prominent candidates like (oriented) homogeneous spaces can be found.

Let us first recall some simple and well-known properties from homological algebra. (By $\zz_p$ we denote the cyclic group of order $p$.)
\begin{lemma}\label{lemma03}
It holds: $\Ext_\zz(A,\qq/\zz)=0$ for all abelian groups $A$ and $\Ext_\zz(\qq,\zz_p)=0$.
\end{lemma}
\begin{prf}
The first assertion is just the injectivity of $\qq/\zz$ as a $\zz$-module. The second assertion follows from the short exact sequence
\begin{align*}
0\to\zz_p\to \qq/\zz \xto{\cdot p} \qq/\zz\to 0
\end{align*}
and its associated sequence
\begin{align*}
0\to \Hom(\qq,\zz_p)&\to \Hom(\qq,\qq/\zz)\xto{\cdot p} \Hom(\qq,\qq/\zz) \\&\to \Ext(\qq,\zz_p)\to \Ext(\qq,\qq/\zz)=0
\end{align*}
However, multiplication with $p$ is surjective on $\Hom(\qq,\qq/\zz)$ and the result follows.
\end{prf}

Let $K$ be a simply-connected finite CW-complex the rational cohomology of which satisfies Poincar\'e Duality together with a rationalisation $K_{(0)}$---which always exists and again is simply-connected. That is, the space $K_{(0)}$ is a rational space, i.e.~$\pi_*(K_{(0)})$ is a $\qq$-vector space, and there is a map $\varphi: K\to K_{(0)}$ which induces an isomorphism
\begin{align*}
\pi_*(X)\otimes \qq\to \pi_*(X_{(0)})
\end{align*}
Suppose that $f: K_{(0)}\to K_{(0)}$ is a continuous endomorphism. We call a map
\begin{align*}
kf: K_{(0)}\to K_{(0)}
\end{align*}
a \emph{$k$-th multiple} of $f$ (for $k \in \nn$) if there is a (finite) homogeneous basis $\{x_i\}$ of $H^{>0}(K_{(0)})$  as a $\qq$-vector space such that the induced map in cohomology satisfies
\begin{align*}
H^*(kf): H^*(K_{(0)},\qq) \to H^*(K_{(0)},\qq)\qquad\qquad H^i(kf)(x_i)=l_{x_i}\cdot k\cdot f(x_i)
\end{align*}
for $l_{x_i}\in \nn$. If $\deg f\neq 0$, the existence of a $k$-th multiple of $f$ for every $k\geq k_0> 0$ evidently shows that $K_{(0)}$ is flexible, in particular. However, in this case even more can be said.
\begin{theo}\label{theo03}
If $f: K_{(0)}\to K_{(0)}$ possesses a $k$-th multiple $kf$ for every $k\geq k_0>0$ and if $\deg f\neq 0$, then $K$ is flexible.
\end{theo}
\begin{prf}
The rational space $K_{(0)}$ is obviously flexible. It remains to see that there exists a $k_0\geq 0$ such that the $k$-th multiple $kf$ of $f$ is covered by a map $\widetilde{kf}: K\to K$ for every $k\geq k_0$.

Since $K_{(0)}$ is the rationalisation of $K$ we obtain a diagram
\begin{align*}
\xymatrix{ K \ar[d]_{h}^\simeq & K \ar[d]_{h}^\simeq\\
K_{(0)}\ar[r]_<<<<{kf} &K_{(0)} }
\end{align*}
with a rational equivalence $h$---cf.~\cite{FHT01}.17.15, p.~256. Hence we aim to construct the lifting $\widetilde{kf}$ of $kf\circ h$.
\begin{align*}
\xymatrix{ & K \ar[d]_{h}^\simeq\\
K\ar[r]_<<<<{kf\circ h} \ar@{.>}[ur]^{\widetilde{kf}}  & K_{(0)}}
\end{align*}
By $F$ we denote the homotopy fibre of the map $h: K\to K_{(0)}$.
Due to obstruction theory there is a sequence of obstructions $c_{kf\circ h}^{n+1}\in H^{n+1}(K, \pi_n(F))$ with the following properties: The obstruction $c_{kf \circ h}^{n+1}$ is defined---depending on a chosen lift---if $c_{kf\circ h}^{n}$ vanishes; and all the obstructions vanish if and only if a lift $\widetilde{k f}$ exists---cf.~theorem \cite{Bre97}.VII.13.11, p.~507.

Let us briefly recall how the obstructions $c_{kf\circ h}^{n}$ are constructed: Consider the Postnikov tower
\begin{align*}
\dots \to K_{n+1} \to K_{n} \to \dots \to K_0=K_{(0)}
\end{align*}
over $K_{(0)}$. It consists of principal fibrations
\begin{align*}
K(\pi_n,n)\hto{} K_{n+1}\xto{h_n} K_{n}
\end{align*}
and identifies $K$ with $\varprojlim K_n$ up to weak equivalence---cf.~theorem \cite{Bre97}.VII.13.7, p.~504, and corollary \cite{Bre97}.VII.13.10, p.~506. By theorem \cite{Bre97}.VII.13.8, p.~505, we have that $\pi_n=\pi_n(F)$.

We shall now successively lift the map $kf\circ h: K\to K_0$ over the subsequent stages in the Postnikov tower. Since every stage is a principal fibration with fibre an Eilenberg-MacLane space $K(\pi_n,n)$, it is classified by a map
\begin{align*}
\theta: K_n \to \B K(\pi_n,n)=K(\pi_n,n+1)
\end{align*}
From \cite{Bre97}, p.~499, respectively from theorem \cite{Bre97}.VII.13.2, p.~500, we see that the obstructions $c_{kf\circ h}^{n}$ are given by
\begin{align*}
c_{kf\circ h}^{n+1}=\big(H^{n+1}(kf\circ h)\circ H^{n+1}(\theta)\big)(u)\in H^{n+1}(K,\pi_{n})
\end{align*}
where $u \in H^{n+1}(K(\pi_n,n+1),\pi_n)$ is a certain class corresponding to an isomorphism $H_n(K(\pi_n,n),\pi_n)\xto{\cong} \pi_n$---cf.~\cite{Bre97}, p.~491.

\vspace{5mm}

In the following we shall show that we may choose $k$ large enough to guarantee that all the obstructions $c_{kf\circ h}^{n}$ vanish successively.

The induced morphism
\begin{align*}
\pi_*(K)\otimes \qq \to \pi_*(K_{(0)}) \otimes \qq=\pi_*(K_{(0)})
\end{align*}
is an isomorphism---cf.~theorem \cite{FHT01}.17.12, p.~254. Since the $\pi_n(K)$ are finitely-generated abelian groups, they split as
\begin{align*}
\pi_n(K)\cong \zz^{l_1}\oplus \zz_2^{l_2} \oplus \zz_3^{l_3}\oplus \zz_5^{l_5}\oplus \dots
\end{align*}
The rationalised long exact homotopy sequence of the fibration
\begin{align*}
F\hto{} K\to K_{(0)}
\end{align*}
shows that $\pi_*(F)$ is only torsion, and so the long exact homotopy sequence yields that
\begin{align*}
\pi_n(F)\cong (\qq/\zz)^{m_1}\oplus \zz_2^{m_2} \oplus \zz_3^{m_3}\oplus \zz_5^{m_5}\oplus \dots
\end{align*}
As cohomology with coefficients a direct sum splits as the direct sum of the cohomology in the respective summands, we hence may simplify the situation. We thus basically may deal with the following two different cases: Either $\pi_n=\qq/\zz$ or $\pi_n=\zz_p$, a finite cyclic group of prime order.

Let us discuss $\pi_n=\zz_p$ first. The obstruction $c_{kf\circ h}^{n+1}$ is the pullback of $u$ along the following morphisms
\begin{align*}
H^{n+1}(K(\pi_n,n+1),\pi_n) &\xto{H^{n+1}(\theta)}  H^{n+1}(K_{(0)},\pi_{n})\\& \xto{H^{n+1}(kf)} H^{n+1}(K_{(0)},\pi_{n}) \\&\xto{H^{n+1}(h)} H^{n+1}(K,\pi_{n})
\end{align*}
Since $H_*(K_{(0)},\zz)$ is a $\qq$-vector space, we compute
\begin{align*}
\Hom(H_i(K_{(0)},\zz_p))=0=\Ext(H_i(K_{(0)},\zz_p)
\end{align*}
using lemma \ref{lemma03} and the fact that both $\Hom$ and $\Ext$ commute with finite sums. Thus the universal coefficient theorem for cohomology (cf.~theorem \cite{Bre97}.V.7.2, p.~282) lets us derive that
\begin{align*}
H^*(K_{(0)},\zz_p)=0
\end{align*}
In particular the pullback of $u$ vanishes.

\vspace{5mm}

Let us now assume that $\pi_n=\qq/\zz$. In this case lemma \ref{lemma03} yields that
\begin{align*}
H^{i}(K,\pi_n)\cong \Hom(H_i(K),\pi_n)\cong (\qq/\zz)^{\tilde m_1} \oplus \zz_2^{\tilde m_2} \oplus \zz_3^{\tilde m_3}\oplus \zz_5^{\tilde m_5}\oplus \dots
\end{align*}
is torsion, as the integral cohomology of a finite CW-complex is a finitely generated abelian group. Thus for every cohomology class $x\in H^{n+1}(K_{(0)},\pi_n)$ there exists a multiple $\tilde k_x \cdot x$ with $\tilde k_x\in \nn$ such that $\tilde k_x \cdot x \in \ker H^{n+1}(h)$.

By assumption the morphism $f$ admits a $k$-th multiple for every $k\geq k_0$. Thus there is a homogeneous basis $\{x_i\}$ of $H^{>0}(K_{(0)})$ such that
\begin{align*}
H^*(kf): H^*(K_{(0)},\qq) \to H^*(K_{(0)},\qq)\qquad\qquad H^i(kf)(x_i)=l_{x_i}\cdot k\cdot f(x_i)
\end{align*}
for $k\geq k_0$.
Applying the observation above to the images under $f$ of the finitely many basis elements $x_i^{n+1}$ in degree $n+1$ we derive that there are coefficients $\tilde k_{x_i^{n+1}}$ with $\tilde k_{x_i^{n+1}} \cdot f(x_i^{n+1}) \in \ker H^{n+1}(h)$. Now set
\begin{align*}
k_{n+1}:=k_0\cdot\prod \tilde k_{x_i^{n+1}}\geq k_0
\end{align*}
Then $k_{n+1} f(x_i^{n+1}) \in \ker H^{n+1}(h)$ and $k_{n+1} f(x) \in \ker H^{n+1}(h)$ for all \linebreak[4]$x\in H^{n+1}(K_{(0)},\pi_n)$. Clearly, the map $k_{n+1}f$ exists and the pullback of $u$ vanishes. Thus there is no obstruction to a lifting in this stage. Now we choose
\begin{align*}
k:=\prod_{0\leq n\leq \dim K-1} k_{n+1}\geq k_0
\end{align*}
This guarantees the existence of $kf$ and, furthermore, it grants that every obstruction to a lifting of $kf$ which successively occurs vanishes. Thus the map $kf$ lifts to $\widetilde{kf}: K\to K$. So does every multiple of $kf$.

(As a final comment on flexibility, we observe that, in particular, for infinitely many $\tilde k\in \nn$ there is the multiple $\tilde kf$ of $f$ which lifts to $\widetilde{\tilde kf}: K\to K$. Every such map $\widetilde{\tilde k f}$ has non-trivial mapping degree, as $\deg f\neq 0$; due to Poincar\'e Duality the mapping degree is well-defined. Thus $K$ is flexible.)
\end{prf}
We give a slightly different view on what we proved, i.e.~actually, a corollary to the proof of theorem \ref{theo03}.
\begin{cor}\label{cor03}
Let $M,N$ be simply-connected finite CW-complexes with rationalisations $M\xto{h_1} M_{(0)}$ and $N\xto{h_2}N_{(0)}$. Let $f_0: M_{(0)} \to N_{(0)}$ be a continuous map. Suppose that $M_{(0)}$ respectively $N_{(0)}$ has a self-map $g$ which admits a $k$-th multiple $g_k$ for each $k\geq k_0>0$. Then there exists a continuous map
\begin{align*}
f&:=\widetilde{f_0\circ g_k \circ h_1}:  M\to N
\intertext{respectively}
f&:=\widetilde{g_k\circ f_0\circ h_1}:  M\to N
\end{align*}
for ($k$ large enough) making the diagram
\begin{align*}
\xymatrix{  M\ar[d]_{h_1}\ar[r]^f & N \ar[d]_{h_2}^\simeq\\
M_{(0)}\ar[r]_<<<<{f_0} &N_{(0)} }
\end{align*}
commute.
\end{cor}
\begin{prf}
Suppose first that $M_{(0)}$ has a self-map $g$ with multiples $g_k$. Adapting the proof of theorem \ref{theo03} we observe that this implies that there is a certain multiple $g_k$ such that the image $\big(H^*(g_k)\circ H^*(f_0)\big)(u)$ of any obstruction class $u$ occurring during the different lifting steps lies in the kernel of $H^*(h_1)$. Thus the obstructions to lifting $h_1\circ g_k\circ f_0$ to a map $f: M\to N$ vanish.

If $N_{(0)}$ has a self-map $g$ with multiples $g_k$, we lift the map $g_k\circ f_0\circ h_1: M\to N_{(0)}$ to a map $f: M\to N$ for $k$ large enough. As before, we may choose $k$ such that the map induced in cohomology $H^*(f_0)\circ H^*(g_k)$ pulls back any obstruction class to an element in the kernel of $H^*(h_1)$.

Clearly, in each case again we draw on the fact that $H^*(M,\pi)$ is a finitely-generated abelian group.
\end{prf}

With theorem \ref{theo03} at hand it is very easy to reprove a key result from \cite{Shi79}.
\begin{cor}
A simply-connected formal space satisfying Poincar\'e Duality is flexible.
\end{cor}
\begin{prf}
We have the following commutative diagram for a formal space $X$ with minimal model $(\bigwedge V,\dif)$:
\begin{align*}
\xymatrix{ (\bigwedge V,\dif) \ar[d]_\simeq  \ar[r]^{\varphi_f}&(\bigwedge V,\dif) \ar[d]^\simeq\\
(H^*(X),0)\ar[r]_<<<<{f} &(H^*(X),0) }
\end{align*}
with the morphism of graded algebras $f$ defined by
\begin{align*}
f: H^k(X)\to H^k(X) \qquad\qquad f(x)=2^{\deg x} x
\end{align*}
and with a Sullivan representative $\varphi_f$ of $f$. The morphism $f$ has a $k$-th multiple $kf$ for every $k\geq 1$ (with respect to an arbitrary homogeneous basis) given by
\begin{align*}
kf(x)=(2k)^{\deg x} x
\end{align*}
So does the realisation $|\varphi_f|: |(\bigwedge V,\dif)|\to |(\bigwedge V,\dif)|$ of $\varphi_f$, by the definition of its $k$-th multiple as
\begin{align*}
k|\varphi_f|:=|\varphi_{kf}|
\end{align*}
and since $H^*(|\varphi_{kf}|)=kf$ without restriction. By theorem \ref{theo03}, the formal space $X$ is flexible, since $\deg f\neq 0$.
\end{prf}

In \cite{CL11}.6.19, p.~20, it was shown that every rational homology class of a pure rational space is the sum of flexible homology classes---cf.~definition \cite{CL11}.6.3, p.~15. We strengthen this in the next theorem.

A Sullivan algebra $(\bigwedge V,\dif)$ (with $V\neq \qq$) is called \emph{two-stage} if it admits a decomposition $V=P\oplus Q$ with $\dif Q=0$ and $\dif P\In \bigwedge Q$. We call a simply-connected topological space \emph{two-stage} if it admits a two-stage Sullivan model. A simply-connected topological space with finite-dimensional rational cohomology admitting a two-stage (minimal) model is necessarily rationally elliptic. Consequently, its rational cohomology satisfies Poincar\'e Duality and there is an orientation class, in particular. So it is legitimate to speak about flexibility in this case.
\begin{theo}\label{theo02}
Let $X$ be a simply-connected two-stage space with finite-dimensional rational cohomology. Then $X$ is flexible.
\end{theo}
\begin{prf}
As above we consider the two-stage decomposition of its minimal model $(\bigwedge V,\dif)$ given by $V=P\oplus Q$ with $\dif Q=0$ and $\dif P\In \bigwedge Q$.
Let $\{q_1,\dots, q_n\}$ be a basis of $Q$ and $\{p_1,\dots, p_m\}$ a basis of $P$.
We now define an automorphism $f$ of graded algebras on $\bigwedge V$ by
\begin{align*}
f(v)&:=2^{\deg v}v
\intertext{for $v\in Q$ and by}
f(v)&:=2^{\deg v+1}v
\end{align*}
for $v\in P$.

As shown in proposition \cite{CL11}.6.19, p.~20, this automorphism is compatible with the differential. Let us now see that this morphism $f$ admits a $k$-th multiple for every $k\geq 1$.

Since $(\bigwedge V,\dif)$ is two-stage, we may decompose it in the form
\begin{align*}
\big(\bigwedge V,\dif\big)&=\bigg(\sum_l \bigwedge Q\otimes \sideset{}{^{\!\! l}}\bigwedge P,\dif\bigg)
\intertext{This yields a decomposition}
H\big(\bigwedge V,\dif\big)&=\sum_i H_i\bigg(\sum_l \bigwedge Q\otimes \sideset{}{^{\!\! l}}\bigwedge P,\dif\bigg)
\end{align*}
in cohomology where an element in $H_i(\sum_l \bigwedge Q\otimes \bigwedge^l P,\dif)$ is represented by a form in $\bigwedge Q\otimes \bigwedge^i P$. (Since the differential preserves the lower grading, i.e.~for $v=\sum v_l$ with $v_l\in \bigwedge Q\otimes \bigwedge^l P$ we have $\dif v=\sum \dif v_l$ with $\dif v_l\in \bigwedge Q\otimes \bigwedge^{l-1} P$, the lower grading of an element a cohomology class is represented by is unique up to exact classes.) Thus we may choose a basis $\{[(x_l^m)_i]\}$ of $H^{>0}(X)$ which is homogeneous in both degrees, i.e.~
\begin{align*}
x_l^m\in \bigg(\big(\bigwedge V\big)_l^m,\dif\bigg)=\bigg(\bigwedge Q\otimes \sideset{}{^{\!\! l}}\bigwedge P\bigg)^m
\end{align*}

Let us now see that $f$ has a $k$-th multiple with respect to this basis for every $k\geq 1$. This follows directly by setting
\begin{align*}
kf(v)&:=(2k)^{\deg v}v
\intertext{for $v\in Q$ and by}
kf(v)&:=(2k)^{\deg v+1}v
\end{align*}
for $v\in P$. (This again defines a morphism of differential graded algebras.) Indeed, on the basis $\{[(x_l^m)_i]\}$ of $H^{>0}(X)$ we compute
\begin{align*}
H^*(kf)([x_l^m]):=(2k)^{l+m}[x_l^m]
\end{align*}
Thus $|kf|$ is a $k$-th multiple of $|f|$ and theorem \ref{theo03} will yield the result once we have seen that $\deg f\neq 0$. This, however, is clear, as the basis elements $[x_l^m]$ are homogeneous in both degrees. So a volume form is represented by a multiple of an element $x_l^m$ already. Such an element is mapped to a non-zero multiple under $f$.
\end{prf}

Let $G$ be a compact connected Lie group and let
$H\In G\times G$ be a closed Lie subgroup. Then $H$ acts on $G$ on
the left by $(h_1,h_2)\cdot g=h_1gh_2^{-1}$. The orbit space of this
action is called the \emph{biquotient} $\biq{G}{H}$ of $G$ by $H$. If the action of $H$ on $G$ is free, then
$\biq{G}{H}$ possesses a manifold structure. This is the only case
we shall consider. Clearly, the category of biquotients contains the
one of homogeneous spaces. It was shown in \cite{Kap}.1, p.~2, that biquotients admit pure models. Trivially, a pure model is two-stage, in particular.
From proposition \cite{Oni94}.1.17, p.~84, we cite that homogeneous
spaces $G/H$ of a compact connected Lie group $G$ with a connected
closed Lie subgroup $H$ are simple.
That is, they have abelian fundamental group acting trivially on higher homotopy groups, and Rational Homotopy Theory and obstruction theory work as well in this setting as in the simply-connected one.
So we have proved
\begin{cor}\label{cor01}
Every simply-connected biquotient is flexible. Besides, so is every (oriented) homogeneous space $G/H$ with $H$ connected.
\end{cor}
\vproof
We remark that there are several simply-connected non-formal homogeneous spaces like $\Sp(n)/\SU(n)$ for $n\geq 5$.

\vspace{5mm}

On a minimal Sullivan algebra $(\bigwedge V,\dif)$ we may consider a lower grading iteratively defined:
As a first step we set $V_0:=\ker \dif|_V$. Suppose that the grading is defined on $V$ until degree $i$, i.e.~we have well-defined vector subspaces $V_1,\dots, V_i$ of $V$ with the property that
\begin{align*}
\dif V_j\in \big(\bigwedge (V_0\oplus \dots \oplus V_{j-1})\big)_{\leq j-1}
\end{align*}
for all $0\leq j\leq i$. Now choose a complement $V_{i+1}$ of
$V_0\oplus \dots \oplus V_{i}$ in $\{v\in V\mid \dif v\in (\bigwedge (V_0\oplus \dots \oplus V_i))_{\leq i}\}$
and extend this grading additively under the multiplicative structure of the algebra to the whole of $(\bigwedge (V_0\oplus \dots \oplus V_{i+1}),\dif)$. Then proceed inductively to finally obtain $V=\oplus_{j\geq 0} V_j$.
\begin{prop}\label{prop04}
Let $M$ be a simply-connected compact closed manifold with minimal model $(\bigwedge V,\dif)$. Suppose that $(\bigwedge V,\dif)$ admits a lower grading as above such that for all $i\in \nn_0$ and for all $v_i\in V_i$ it holds that
\begin{align*}
\dif(v_i)\in \big(\bigwedge (V_0\oplus \dots \oplus V_{i-1})\big)_{i-1}=\big(\bigwedge V\big)_{i-1}
\end{align*}
Then $M$ is flexible.
\end{prop}
\begin{prf}
We define a morphism of commutative graded algebras $f: (\bigwedge V,\dif) \to (\bigwedge V,\dif)$ on $V_{i}^j$ for $i,j\geq 0$ by
\begin{align*}
f(v):=2^{i+j} \cdot v
\end{align*}
The multiplicative extension of the grading implies that for $v\in (\bigwedge V)_i^j$ we then have
\begin{align*}
f(v)=2^{i+j}\cdot v
\end{align*}
This morphism is obviously compatible with the differential. Indeed, for an arbitrary element in $(\bigwedge V,\dif)$ given by
\begin{align*}
v&=\sum_m (v_m)_i^j \textrm{\qquad\qquad with \qquad\qquad} (v_m)_i^j\in \big(\bigwedge V\big)_i^j
\intertext{we compute}
\dif(f(v))&=\dif\bigg(\sum_m 2^{i+j} \cdot (v_m)_i^j\bigg)=\bigg(\sum_i 2^{(i-1)+(j+1)} \cdot \dif \big((v_m)_i^j\big)\bigg)=f(\dif v)
\end{align*}
since $\dif (v_m)_i^j\in (\bigwedge V)_{i-1}^{j+1}$.

The double grading on $(\bigwedge V,\dif)$ induces  a double grading on cohomology, i.e.~we have a decomposition
\begin{align*}
H\big(\bigwedge V,\dif\big)&=\sum_{i,j} H_i^j\big(\bigwedge V,\dif\big)
\intertext{where}
H_i^j\big(\bigwedge V,\dif\big)
\end{align*}
is the subgroup of $H^j(\bigwedge V,\dif)$ represented by closed elements in $(\bigwedge V)_i$. We remark that this is well-defined, since the differential decreases the lower degree by exactly one by assumption.

Now choose a basis of $H(\bigwedge V,\dif)$ subordinate to the bigrading on cohomology. The same arguments as for $f$ above show that the morphism $f$ obviously has a $k$-th multiple $kf$ for $k\geq 1$ well-defined by
\begin{align*}
kf(v)=(2k)^{i+j}\cdot v
\end{align*}
for $v\in (\bigwedge V)_i^j$ with respect to such a basis. Due to the bigrading on cohomology the mapping degree of $f$ is non-zero and $M$ is flexible by theorem \ref{theo03}.
\end{prf}

\begin{rem}
The \emph{bigraded model} of a formal space (cf.~theorem \cite{FOT08}.2.93, p.~95) has the properties required in proposition \ref{prop04}. Consequently, this reproves that simply-connected formal spaces are flexible. However, a space with such a special grading needs not be formal as we shall see in example \ref{ex02}.

Besides, the we see that the grading we used on two-stage spaces in theorem \ref{theo02} has all the properties of proposition \ref{prop04}, Thus this proposition also comprises the result for two-stage spaces. However, for the convenience of the reader we provided the respective special proofs.

Compare our lower grading to the setting of \emph{universal spaces}---cf.~\cite{AL00} and lemma \cite{AL00}.4.5, p.~534---in which everything works equally well.
\end{rem}

\begin{cor}\label{cor02}
Let $M$ be a simply-connected compact closed manifold with minimal Sullivan model $(\bigwedge V,\dif)$.
Suppose that there is a basis $\{v_i\}$ of $V$ such that for each $i$ the differential $\dif v_i$ of a basis element is a monomial in some other basis elements $v_j$. Then $M$ is flexible.
\end{cor}
\begin{prf}
By a change of basis of $V$, which does not affect the fact that differentials are monomials, we may assume that the given basis is subordinate to $\ker \dif|_V$, i.e.~there are basis elements $v_j$ which span $\ker \dif|_V$.

We define a lower grading on $(\bigwedge V,\dif)$ as above---however, we make sure that the basis given in the assertion is subordinate to the grading. Thus the complements are not chosen arbitrarily: So we choose for $V_0$ the space spanned by all the $v_i$ with $\dif v_i=0$. We set $V_1$ to be the complement of $V_0$ in $V$ spanned by all the $v_i$ which do not span $V_0$ and which map to $\bigwedge V_0$. Iteratively, we have $V_{i+1}$ be spanned by the $v_i$ which do not form an element of any $V_j$ with $j\leq i$ and which map to $(\bigwedge (V,\dif))_{\leq i}$.

It is clear that the basis elements by which we span the $(i+1)$-th degree $V_{i+1}$ really form a complement of $V_\leq i$ in $\{v\in V\mid \dif v\in (\bigwedge (V_0\oplus \dots \oplus V_i))_{\leq i}\}$, i.e.~there is no non-trivial linear combination of them which has lower degree. This is due to the fact that the $\dif v_i$ are monomials in the basis and that we displayed $\ker \dif|_V$ with maximal dimension in the basis elements. By the same arguments we see that actually
\begin{align*}
\dif(v_i)\in \big(\bigwedge (V_0\oplus \dots \oplus V_{i-1})\big)_{i-1}=\big(\bigwedge V\big)_{i-1}
\end{align*}
not only for all basis elements $v_i\in V_i$ but for all elements $v_i\in V_i$. By proposition \ref{prop04} the result follows.
\end{prf}

\begin{ex}\label{ex02}
The examples we used to prove theorem \ref{theoA} were three-stage spaces, i.e.~their minimal model $(\bigwedge V,\dif)$ could be decomposed as $V=V_1\oplus V_2\oplus V_3$ with $\dif V_i\in \bigwedge V_{\leq i}$ or, equivalently, their rational Postnikov tower has height $3$. As we have proved in theorem \ref{theo02} a two-stage space is flexible although not necessarily formal. So let us provide an easy example of a flexible three-stage simply-connected smooth compact closed manifold which is not formal. This will be an easy application of corollary \ref{cor02}.

For this we define a minimal Sullivan algebra $(\bigwedge V,\dif)$ by
\begin{align*}
V=\langle a,b,n,m\rangle
\end{align*}
with $\deg a=\deg b=3$, $\deg n=5$, $\deg m=7$. Define the differential by $\dif a=\dif b=0$, $\dif n=ab$, $\dif m=an$. Since all the finitely many generators are of odd degree, this defines an elliptic algebra. Its dimension is $3+3+5+7=18$. Consequently, it is realisable by a simply-connected smooth compact closed manifold $M$. The manifold $M$ is obviously not formal, since the class $bn$ is closed but not exact. The basis $\{a,b,n,m\}$ has all the properties required in corollary \ref{cor02}; thus, the manifold $M$ is flexible.

In order to illustrate the used arguments again, we give a short and direct proof in the spirit of proposition \ref{prop04}: Choose the basis
\begin{align*}
[a], [b], [bn], [am], [anm], [bnm], [abmn]
\end{align*}
for its cohomology of positive degree, we may define the following morphism of differential graded algebras on $(\bigwedge V,\dif)$ which has multiples and therefore yields that $M$ is flexible.
\begin{align*}
f(a)&:= 2^3 a \\
f(b)&:= 2^3 b \\
f(n)&:= 2^6 n \\
f(m)&:= 2^9 m \\
\end{align*}
A $k$-th multiple of $f$ for $k\geq 1$ with respect to the basis above is easily seen to be
\begin{align*}
kf(a)&:= (2k)^3 a \\
kf(b)&:= (2k)^3 b \\
kf(n)&:= (2k)^6 n \\
kf(m)&:= (2k)^9 m \\
\end{align*}
and by theorem \ref{theo03} the manifold $M$ is flexible.
\end{ex}

\vspace{5mm}
In \cite{AL00} a conjecture of Copeland--Shar saying that for rational spaces $M$ and $N$ the set $[M,N]$ of homotopy classes of maps is either infinite or consists of a single element only. However, for a large class of spaces, namely \emph{universal spaces}, the conjecture is proven to hold true in theorem \cite{AL00}.4.6, p.~535. See definition \cite{AL00}.4.2, p.~534, for the precise definition of universality. We shall just mention that this class comprises formal spaces, Eilenberg--MacLane spaces, Moore spaces, H-spaces and, in particular, homogeneous spaces besides many more. We shall no ``lift'' this conjecture to the ``real world''.
\begin{theo}
Let $M$, $N$ be simply-connected finite CW-complexes. Suppose that either $M$ or $N$ is universal and that $[M_{(0)},N_{(0)}]\neq [\id]$ contains a non-trivial element. Then $\big|[M,N]\big|=\infty$.
\end{theo}
\begin{prf}
Since a finite complex is universal if and only if so is its minimal model---cf.~\cite{AL00}, p.~534; or, equivalently, its rationalisation.
Denote by $(\bigwedge V,\dif)$ the minimal model of $M$ and by $(\bigwedge W,\dif)$ the one of $N$.

Suppose first that $M$ is universal.
In theorem \cite{AL00}.4.6, p.~535, it is proved that $[M_{(0)},N_{(0)}]$ contains infinitely many elements if it contains more than one. Indeed, this is done by composing the self-maps $\phi_k$ defined by $\phi_k(v)=k^{\deg v+i}v$ for $v\in V_i$---which exist due to universality---with a certain non-trivial map $f_0: M_{(0)}\to N_{(0)}$. Obviously, a $(k^{i-1})$-th multiple of $\phi_k$ is given by $(\phi_k)^i$. (We use a cohomology basis compatible with the used bidegree---cf.\cite{AL00}, p.~534 and lemma \cite{AL00}.4.5 on that page.) In the terminology of corollary \ref{cor03} we thus obtain a lift of $f_0\circ (\phi_k)^i \circ h_1: M\to N_{(0)}$ for all $i\geq i_0$ for a certain large $i_0\in \nn$; i.e.~for infinitely many such maps. As was proved in \cite{AL00}, the maps $f_0\circ \phi_k$ are pairwise distinct elements of $[M_{(0)},N_{(0)}]$. So the maps $f=\widetilde{f_0\circ (\phi_k)^i \circ h_1}$ will represent pairwise distinct classes in $[M,N]$. (Indeed, a homotopy of maps $f$ induces a homotopy of Sullivan representatives $f_0\circ \phi_k$---cf.~proposition \cite{FHT01}.12.6, p.~149. The realisations of two such representatives are homotopic if so are the representatives themselves---cf.~proposition \cite{FHT01}.17.13, p.~255.). This proves that $[M,N]$ has infinitely many elements.

In the case when $N$ is universal, we apply analogous arguments to the compositions $(\phi_k)^i \circ f_0\circ h_1$. Corollary \ref{cor03} and the fact that the $(\phi_k)^i\circ f_0$ are pairwise not homotopic according to the proof of theorem \cite{AL00}.4.6, p.~535, again yield the result.
\end{prf}

%%%%%%%%%%%%%%%%%%%%%%%%%%%%%%%%%% Section 4 %%%%%%%%%%%%%%%%%%%%%%%%%%%%%%%%%%%%%%

\section{Manifolds without orientation reversal}\label{sec04}

In \cite{Mue09}, theorem B, p.~2362, it was shown that from dimension $7$ onwards in each dimension a simply-connected \emph{strongly chiral} smooth manifold, i.e.~a manifold which does not possess a self-map of degree $-1$, exists. In dimensions divisible by four the complex projective spaces yield very simple examples. In \cite{Mue09}, p.~2370, it is illustrated that in each dimension congruent to $3$ modulo $4$---starting in dimension $7$---there is an $\s^{2k-1}$-bundle over $\s^{2k}$ the total space of which does not permit an orientation reversal. However, in the other cases the construction of simply-connected strongly chiral manifolds relies heavily on Cartesian products.

We shall generalise these results in two directions: As announced we shall consider arbitrary negative mapping degree and we shall search for irreducible simply-connected examples. We provide the following result which shows that from a certain dimension on irreducible examples can be found in infinitely many dimensions.
Moreover, we remark that the following examples lack orientation reversal although they are not inflexible. This is due to the fact that their minimal models are pure algebras---cf.~theorem \ref{theo02}. Let us first deal with the case of $\dim M\equiv 2 \mod 4$.
\begin{prop}\label{prop01}
In every dimension $4k+54\geq 58$ ($k\in \nn$) there is an irreducible simply-connected smooth compact closed manifold which does not admit any self-map of negative degree.
\end{prop}
\begin{prf}
Define a minimal Sullivan algebra $(\bigwedge V,\dif)$ by
\begin{align*}
V:=\langle x_1,x_2,n_1,n_2,n_3,n_4\rangle
\end{align*}
with $\deg x_1=2, \deg x_2=4, \deg n_1=11, \deg n_2=11, \deg n_3=4l_1+1,\deg n_4=8 l_2+3$ and
\begin{align*}
\dif x_1&=\dif x_2=0\\
\dif n_1&=x_1^4x_2 \\
\dif n_2&=x_1^2x_2^2+x_2^3\\
\dif n_3&=x_1^{2l_1+1}\\
\dif n_4&=x_2^{2l_2+1}
\end{align*}
with $l_1, l_2\geq 2$. This algebra is finitely-generated. Since $x_1,x_2$ define nilpotent cohomology classes, the algebra is obviously has finite-dimensional cohomology. Thus it is elliptic and satisfies Poincar\'e duality. Hence we may compute its dimension (cf.~\cite{FHT01}.32, p.~434) as
\begin{align*}
11+11+4l_1+1+8l_2+3-1-3=4l_1+8l_2+22\geq 46
\end{align*}
Thus $\dim (\bigwedge V,\dif)\equiv 2\mod 4$. Consequently, this algebra may be realised as the minimal model of a simply-connected smooth compact closed manifold $M$. As the model does not split non-trivially as a product, neither does the manifold.

\vspace{5mm}

A volume form is given by
\begin{align*}
v:=x_1^{2l_1}x_2^{2l_2}n_1n_2- (x_1^{2l_1+2}x_2+x_1^{2l_1}x_2^{2}) n_1n_4 + x_1^{2l_1+4}n_2n_4
\end{align*}
Indeed, we compute
\begin{align*}
&\dif\big(x_1^{2l_1}x_2^{2l_2}n_1n_2- (x_1^{2l_1+2}x_2+x_1^{2l_1}x_2^{2}) n_1n_4 + x_1^{2l_1+4}n_2n_4\big)\\
=&x_1^{2l_1+4}x_2^{2l_2+1}n_2-x_1^{2l_1+2}x_2^{2l_2+2}n_1-x_1^{2l_1}x_2^{2l_2+3}n_1
\\&-x_1^{2l_1+6}x_2^{2}n_4-x_1^{2l_1+4}x_2^{3}n_4+x_1^{2k+2}x_2^{2l_2+2}n_1+x_1^{2l_1}x_2^{2l_2+3}n_1\\
&+x_1^{2l_1+6}x_2^2n_4+x_1^{2l_1+4}x_2^3n_4-x_1^{2k_1+4}x_2^{2k_2+1}n_2
\\=&0
\end{align*}
The class cannot be exact, as every class the differential of which contains $x_1^{2l_1}x_2^{2l_2}n_1n_2$ as a summand has to contain either a summand with factors $n_1n_2n_3$ or one which has factors $n_1n_2n_4$. In the first case differentiating yields a power too large in $x_1$, in the second case a power too large in $x_2$.

\vspace{5mm}

We shall now consider an arbitrary morphism $f$ on this model. It is of the form
\begin{align*}
f(x_1)&=k_1x_1\\
f(x_2)&=k_2 x_1^2+k_3 x_2\\
f(n_1)&=k_4 n_1+k_5 n_2\\
f(n_2)&=k_6 n_1+k_7 n_2
\intertext{Since $f$ commutes with the differential, it follows that}
f(n_3)&=k_8n_3+c\\
f(n_4)&=k_9n_4+d
\end{align*}
with $\dif c=\dif d=0$.

Indeed, the algebra $(\bigwedge V,\dif)$ is pure. Thus we may consider the lower grading which is respected by the differential. In particular, this shows that the differential $\dif x$ of an element $x\in \bigwedge V$ may only have a power in $x_1$ or in $x_2$ as a summand if $v\in \bigwedge_1 V$, i.e.~every summand of $v$ has at most one factor being an $n_i$. However, $p\dif n_i=x_i^l$ (for $l>0$) with $p\in \bigwedge\langle x_1,x_2\rangle$ implies that $n_i=n_3$ respectively $n_i=n_4$. From here we obtain the above structure of $f$ on $n_3$ and $n_4$.

We shall now focus on the case when $l_2=2$ and $l_1\geq 4$. This implies that $d=0$ by degree, as a little check directly shows.
The formal dimension of the algebra then is
\begin{align*}
\dim(A,\dif)=38+4l_1\geq 54
\end{align*}

Besides, we may compute that
\begin{align*}
f(\dif n_1)&=f(x_1^4x_2)=k_1^4k_2x_1^6+k_1^4k_3x_1^4x_2
\intertext{equals}
\dif(f (n_1))&=\dif(k_4 n_1+k_5 n_2)=k_4x_1^4x_2+k_5x_1^2x_2^2+k_5x_2^3
\end{align*}
which yields $k_4=k_1^4k_3$, $k_5=0$ and $k_1=0\vee k_2=0$. If $k_1=0$ then the volume form $v$ is mapped to zero by $f$ and $f$ does not reverse orientation. Thus we may assume that $k_1\neq 0$ and $k_2=0$ and $f(x_2)=k_3x_2$. Analogously, it follows that
\begin{align*}
f(\dif n_2)&=f(x_1^2x_2^2+x_2^3)=k_1^2k_3^2x_1^2x_2^2+k_3^3x_2^3
\intertext{equals}
\dif(f (n_2))&=\dif(k_6 n_1+k_7 n_2)=k_6x_1^4x_2+k_7x_1^2x_2^2+k_7x_2^3
\end{align*}
We infer that $k_7=k_3^3=k_1^2k_3^2$ and $k_6=0$. So either $k_3=k_4=k_7=0$ or $k_3=k_1^2, k_7=k_1^6$ and $k_4=k_1^6$. Moreover, we see that
\begin{align*}
\dif (f(n_4))&=k_9n_4
\intertext{equals}
f(\dif(n_4))&=f(x_2^5)=k_3^5x_2^5
\end{align*}
which gives us $k_9=k_3^5=k_1^{10}$.

\vspace{5mm}

Summing up these results, we either have $k_3=k_4=k_7=0$ or $k_3=k_1^2$, $k_4=k_7=k_1^6$, $k_9=k_1^{10}$. The first case yields $f(v)=0$ and in the second case we may compute the image of the volume form as
\begin{align*}
f(v)=&f\big(x_1^{2l_1}x_2^{4}n_1n_2- (x_1^{2l_1+2}x_2+x_1^{2l_1}x_2^{2}) n_1n_4 + x_1^{2l_1+4}n_2n_4\big)\\
=&k_1^{2l_1+20}(x_1^{2l_1}x_2^4n_1n_2-(x_1^{2l_1+2}x_2+x_1^{2l_1}x_2^2)n_1n_4+x_1^{2l_1+4}n_2n_4)\\
=&k_1^{2l_1+20} v
\end{align*}
Since $2l_1+20$ is even, the coefficient $k_1^{2l_1+20}$ is non-negative. Thus $f$ cannot reverse orientation. As $f$ was chosen arbitrarily, its realisation as a manifold, $M$, does not admit orientation reversal.
\end{prf}

Now we construct examples in dimensions $\dim M\equiv 1 \mod 4$.
\begin{prop}\label{prop02}
In every dimension $4k+69\geq 73$ ($k\in \nn$) there is an irreducible simply-connected smooth compact closed manifold which does not admit any map of negative degree.
\end{prop}
\begin{prf}
Define a minimal Sullivan algebra $(\bigwedge V,\dif)$ by
\begin{align*}
V:=\langle x_1,x_2,n_1,n_2,n_3,n_4, n_5\rangle
\end{align*}
with $\deg x_1=2, \deg x_2=4, \deg n_1=13, \deg n_2=11, \deg n_3=4l+1,\deg n_4=19, \deg n_5=17$ and
\begin{align*}
\dif x_1&=\dif x_2=0\\
\dif n_1&=x_1^5x_2 \\
\dif n_2&=x_1^2x_2^2+x_2^3\\
\dif n_3&=x_1^{2l+1}\\
\dif n_4&=x_2^{5}\\
\dif n_5&=x_1^3x_2^3
\end{align*}
with $l\geq 4$. Again, this algebra is finitely-generated. Since $x_1,x_2$ define nilpotent cohomology classes, the algebra is obviously has finite-dimensional cohomology. Thus it is elliptic and satisfies Poincar\'e duality. Its dimension computes as
\begin{align*}
13+11+4l+1+19+17-1-3=4l+57\geq 73
\end{align*}
Thus $\dim (\bigwedge V,\dif)\equiv 1\mod 4$. Consequently, also this algebra may be realised as the minimal model of a simply-connected smooth compact closed manifold $M$. As the model does not split non-trivially as a product, neither does the manifold.

\vspace{5mm}

A volume form is given by
\begin{align*}
v:=&x_1^{2l}x_2^4n_1n_2n_5-x_1^{2l+3}x_2^2n_1n_2n_4-(x_1^{2l+2}x_2+x_1^{2l}x_2^2)n_1n_4n_5\\&+x_1^{2l+5}n_2n_4n_5
\end{align*}
Indeed, a direct computation again shows that $v$ is a closed form. Once more, the class cannot be exact, as every class the differential of which contains $x_1^{2l_1}x_2^4n_1n_2n_5$ as a summand has to contain either a summand with factors $n_1n_2n_3n_5$ or one which has factors $n_1n_2n_4n_5$. In the first case differentiating yields a power too large in $x_1$, in the second case a power too large in $x_2$.

\vspace{5mm}

As in the proof of proposition \ref{prop01} we shall now consider an arbitrary morphism $f$ on this model. However, using that it commutes with the differential shows that it is actually of the form
\begin{align*}
f(x_1)&=k_1x_1\\
f(x_2)&=k_2 x_2\\
f(n_1)&=k_3 n_1\\
f(n_2)&=k_4 n_2\\
f(n_3)&=k_5n_3+c\\
f(n_4)&=k_6n_4\\
f(n_5)&=k_7 n_5
\end{align*}
with $c\in \bigwedge V$. Indeed, as above we see that $f(\dif n_1)=\dif(f(n_1)$ implies that $f(n_1)=k_3n_1$ and $f(x_2)=k_2x_2$. Moreover, it is easy to check that there are no non-trivial closed combinations of degree $19=\max \deg_{i\neq 3} n_i$ or smaller of the $\dif n_i$ ($1\leq i\leq 5$, $l\geq 4$) with coefficients in $\bigwedge V$---by degree actually in $\qq[x_1,x_2]$. This shows that each of the $n_i$ ($i\neq 3$) has to map to a scalar multiple of itself.

The fact that $f$ commutes with the differential then gives us $k_3=k_1^5k_2$, $k_4=k_1^2k_2^2=k_2^3$ and $k_2=k_1^2 \vee k_4=k_2=0$, $k_6=k_2^5$, $k_7=k_1^3k_2^3$. (Again, the case $k_2=k_3=k_4=k_6=k_7=0$ cannot occur unless the mapping degree of $f$ is $0$.)

Thus we may compute the image of the volume form as
\begin{align*}
&f(v)\\=&f\big(x_1^{2l}x_2^4n_1n_2n_5-x_1^{2l+3}x_2^2n_1n_2n_4-(x_1^{2l+2}x_2+x_1^{2l}x_2^2)n_1n_4n_5\\&+x_1^{2l+5}n_2n_4n_5\big)
\\=& k_1^{2l+30}v
\end{align*}
Since $2l_1+30$ is even, the coefficient $k_1^{2l_1+30}$ is non-negative. Thus $f$ cannot reverse orientation. As $f$ was chosen arbitrarily, its realisation as a manifold, $M$, does not admit orientation reversal.
\end{prf}
We remark that the minimal models of the manifolds we constructed are pure.

The next examples will satisfy $\dim M\equiv 3 \mod 4$. We shall merely sketch details whenever arguments are similar to what we did above.
\begin{prop}\label{prop03}
In every dimension $4k+43\geq 47$ ($k\in \nn$) there is an irreducible simply-connected smooth compact closed manifold which does not admit any map of negative degree.
\end{prop}
\begin{prf}
The elliptic minimal Sullivan algebra given by \linebreak[4] $V=\langle x_1,x_2,n_1,n_2,n_3\rangle$ with $\deg x_1=2, \deg x_2=4, \deg n_1=11, \deg n_2=19, \deg n_3=4l+1$ and with
\begin{align*}
\dif x_1=\dif x_2&=0\\
\dif n_1&=x_1^2x_2^2+x_2^3\\
\dif n_2&=x_2^5\\
\dif n_3&=x_1^{2l+1}
\end{align*}
has dimension $4l+27\equiv 3 \mod 4$. For $l\geq 5$ (and dimension greater than or equal to $47$ its realisation will yield a manifold without orientation reversal. For this we compute a volume form as
\begin{align}\label{eqn05}
x_2^4x_1^{2l}n_1-x_1^{2l+2}x_2n_2-x_1^{2l}x_2^2n_2
\end{align}
Indeed, this form is closed and non-exact. For the latter we
observe that by the lower grading every exact form of word-length one in the $n_i$ must lie in the $\qq[x_1,x_2]$ span of the $n_1n_2$, $n_1n_3$, $n_2n_3$ with
\begin{align*}
\dif(n_1n_2)&=(x_1^2x_2^2+x_2^3)n_2-x_2^5n_1\\
\dif(n_1n_3)&=(x_1^2x_2^2+x_2^3)n_3-x_1^{2l+1}n_1\\
\dif(n_2n_3)&=x_2^5n_3-n_2x_1^{2l+1}
\end{align*}
Let $v$ be such an element. Consequently, every summand of $v$ of the form $x_1^mn_i$ for some $m>0$ satisfies $m\geq 2l+1$ as a direct computation shows. Hence for $v$ the form from \eqref{eqn05}, we see that it really defines a volume form.

Assume that $l\geq 5$. Every endomorphism $f$ of the model has the form
\begin{align*}
f(x_1)&=k_1 x_1\\
f(x_2)&=k_2 x_2+k_3x_1^2\\
f(n_1)&=k_4n_1\\
f(n_2)&=k_5n_2+k_6x_1^4n_1+k_7x_1^2x_2n_1+k_8x_2^2n_1
\end{align*}
The fact that $f\circ \dif=\dif \circ f$ applied to $n_1$ gives us the following relations:
\begin{align}
\label{eqn06} k_4&=k_1^2k_2^2+3k_2^2k_3\\
\label{eqn07} k_4&=k_2^3\\
\label{eqn08} 2k_1^2k_2k_3+3k_2k_3^2&=0\\
\label{eqn09} k_1^2k_3^2+k_3^3&=0
\end{align}
Combining \eqref{eqn06} with \eqref{eqn07} we obtain that $k_2=k_1^2+3k_3$ unless $k_2=0$. Equation \eqref{eqn09} shows that $k_3=-k_1^2 \vee k_3=0$.

Assume that $k_3\neq 0$.
If $k_2=k_1^2+3k_3$, we derive the implication
\begin{align*}
k_3=-k_1^2\Rightarrow k_1=0
\end{align*}
from equation \eqref{eqn08} and the mapping degree of $f$ is zero.
If $k_2=0$, we obtain $k_3=-k_1^2$, $k_4=0$ and the image of the volume form is
\begin{align*}
&f(x_2^4x_1^{2l}n_1-x_1^{2l+2}x_2n_2-x_1^{2l}x_2^2n_2)
\\=&0+k_1^{2l+4}x_1^{2l+4}f(n_2)-k_1^{2l+4}x_1^{2l+4}f(n_2)
\\=&0
\end{align*}
Altogether, we derive that $k_3=0$ unless the mapping degree of $f$ vanishes. From equations \eqref{eqn06} and \eqref{eqn07} we see that $k_3=0$ implies $k_2=k_1^2$ and $k_4=k_1^6$ unless $k_2=0$, which again yields $\deg f=0$.

An easy calculation using the fact that the differential commutes with $f$ shows that $k_6=k_7=k_8=0$ and that $k_5=k_2^5=k_1^{10}$ unless $\deg f=0$. In this case we thus have
\begin{align*}
k_2&= k_1^2\\
k_3=k_6=k_7=k_8&=0\\
k_4&=k_1^6\\
k_5&=k_1^{10}
\end{align*}
We then may compute the image of the volume form as
\begin{align*}
f\big(x_2^4x_1^{2l}n_1-x_1^{2l+2}x_2n_2-x_1^{2l}x_2^2n_2\big)
=k_1^{2l+14}\big(x_2^4x_1^{2l}n_1-x_1^{2l+2}x_2n_2-x_1^{2l}x_2^2n_2\big)
\end{align*}
Again, there is only an even power of $k_1$ involved and $k_1^{2l+14}$ is positive. Thus the corresponding manifold does not admit any orientation reversal.
\end{prf}

The next triviality provides examples in dimensions $\dim M\equiv \mod 4$.
\begin{rem}\label{rem01}
The complex projective spaces $\cc\pp^{2n}$ form examples which do not permit self-maps of negative degree, since their orientation class is an even power of the K\"ahler class.
\end{rem}

\begin{lemma}\label{lemma02}
There are simply-connected irreducible manifolds of the same rational homotopy type as the manifolds of propositions \ref{prop01}, \ref{prop02} and \ref{prop03} as well as theorem \ref{theo01} which do not decompose non-trivially as connected sums (of oriented compact closed manifolds).
\end{lemma}
\begin{prf}
Decompose a manifold $M$ from the theorems as $M=M'\#N$. The fundamental group of $M$ is the free product of the fundamental groups of $M'$ and $N$. Thus the simple-connectedness of $M$ implies the one of both factors.

We shall now show that the structure of the minimal models of the examples shows that this decomposition has to be rationally trivial, i.e.~$N$ has the rational homotopy type of a sphere: We may assume that $H(M,\qq)$ has a non-trivial cohomology class---other than a multiple of the orientation class being understood. Now suppose that also $H(N,\qq)$ possesses a non-trivial class. A model of $M\#N$ is given by the direct sum of the minimal models of $M$ and $N$ with units and orientation classes identified---cf.~\cite{FHT01}, p.~137---i.e.~by
\begin{align*}
\bigg(\big(\bigwedge V_{M'}\oplus_\qq \bigwedge V_{N'}\big)\oplus \langle u\rangle,\dif\bigg)
\end{align*}
with $\dif(u)=\omega_M-\omega_N$ (where $\omega_M, \omega_N$ are orientation forms) and with $\dif$ restricting to the differential of the respective minimal models on each summand.

The existence of a non-trivial rational cohomology class of $N$ implies that either $x_1$ or $x_2$ lies in $V_N$ in each respective example: This is due to the fact that $x_1,x_2\in V_{M'}$ implies that $n_i\in V_{M'}$ for all $i$ in the case of propositions \ref{prop01}, \ref{prop02} and \ref{prop03}; respectively that all the $y_i,z,z'\in V_{M'}$ in the case of theorem \ref{theo01}. Moreover, it shows that $(\bigwedge V_M,\dif)=(\bigwedge V_{M'},\dif)$ respectively $H(N,\qq)=H\big(\s^{\dim M}\big)$ in each case; a contradiction.

As $H(M,\qq)$ was assumed to be non-trivial, we obtain that, without restriction, $[x_1]\in H(M,\qq)$ and $[x_2]\in H(N,\qq)$. This forces $[x_1x_2]=0$, which is not the case; in neither example. Thus we obtain $H(N,\qq)=H(\s^{\dim M},\qq)$.

This, however, implies that the collapse map $M'\#N\to M'$ induces a quasi-isomorphism. Thus $M'$ has the rational homotopy type of $M$.
\end{prf}
\begin{rem}\label{rem02}
The rational homotopy type of the manifold $(M')^n$ in lemma \ref{lemma02} is the same as the one of $M^n$. All our arguments in the proofs of theorem \ref{theo01} and the propositions \ref{prop01}, \ref{prop02} and \ref{prop03}, in lemma \ref{lemma01} and in remark \ref{rem01} only use the rational homotopy information. So they apply as well to the manifolds $(M')^n$. This implies, in particular, that properties like inflexibility or lack of orientation reversal are passed on from $M$ to $M'$.
\end{rem}
This enables us to give the
\begin{proof}[\textsc{Proof of theorem \ref{theoB}}]
 In all congruence classes modulo $4$ of dimensions starting with dimension $70$ we have constructed irreducible simply-connected compact closed smooth manifolds $M^n$ without orientation reversal in propositions \ref{prop01}, \ref{prop02} and \ref{prop03} respectively in remark \ref{rem01}. We make use of lemma \ref{lemma02} to possibly decompose these examples into simply-connected manifolds $(M')^n$ that do no longer split as connected sums. Due to remark \ref{rem02} the manifolds $(M')^n$ do not possess any self-map of negative degree.
\end{proof}

%%%%%%%%%%%%%%%%%%%%%%%%%%%%%%%%%% Bibliography %%%%%%%%%%%%%%%%%%%%%%%%%%%%%%%%%%%

%\bibliography{lib}{}

\begin{thebibliography}{10}

\bibitem{AL00}
M.~Arkowitz and G.~Lupton.
\newblock Rational obstruction theory and rational homotopy sets.
\newblock {\em Math. Z.}, 235(3):525--539, 2000.

\bibitem{Bre97}
G.~E. Bredon.
\newblock {\em Topology and geometry}, volume 139 of {\em Graduate Texts in
  Mathematics}.
\newblock Springer-Verlag, New York, 1997.
\newblock Corrected third printing of the 1993 original.

\bibitem{CV11}
C.~Costoya and A.~Viruel.
\newblock Every finite group is the group of self homotopy equivalences of an
  elliptic space.
\newblock arXiv:1106.1087v1, 2011.

\bibitem{CL11}
D.~Crowley and C.~Loeh.
\newblock Functorial semi-norms on singular homology and (in)flexible
  manifolds.
\newblock arXiv:1103.4139v1, 2011.

\bibitem{FHT91}
Y.~F{\'e}lix, S.~Halperin, and J.-C. Thomas.
\newblock Elliptic spaces.
\newblock {\em Bull. Amer. Math. Soc. (N.S.)}, 25(1):69--73, 1991.

\bibitem{FHT01}
Y.~F{\'e}lix, S.~Halperin, and J.-C. Thomas.
\newblock {\em Rational homotopy theory}, volume 205 of {\em Graduate Texts in
  Mathematics}.
\newblock Springer-Verlag, New York, 2001.

\bibitem{FOT08}
Y.~F{\'e}lix, J.~Oprea, and D.~Tanr{\'e}.
\newblock {\em Algebraic models in geometry}, volume~17 of {\em Oxford Graduate
  Texts in Mathematics}.
\newblock Oxford University Press, Oxford, 2008.

\bibitem{Kap}
V.~Kapovitch.
\newblock A note on rational homotopy of biquotients.
\newblock preprint.

\bibitem{Mue09}
D.~M{\"u}llner.
\newblock Orientation reversal of manifolds.
\newblock {\em Algebr. Geom. Topol.}, 9(4):2361--2390, 2009.

\bibitem{Oni94}
A.~L. Onishchik.
\newblock {\em Topology of transitive transformation groups}.
\newblock Johann Ambrosius Barth Verlag GmbH, Leipzig, 1994.

\bibitem{Pup95}
V.~Puppe.
\newblock Simply connected {$6$}-dimensional manifolds with little symmetry and
  algebras with small tangent space.
\newblock In {\em Prospects in topology ({P}rinceton, {NJ}, 1994)}, volume 138
  of {\em Ann. of Math. Stud.}, pages 283--302. Princeton Univ. Press,
  Princeton, NJ, 1995.

\bibitem{Shi79}
H.~Shiga.
\newblock Rational homotopy type and self-maps.
\newblock {\em J. Math. Soc. Japan}, 31(3):427--434, 1979.

\bibitem{Sul77}
D.~Sullivan.
\newblock Infinitesimal computations in topology.
\newblock {\em Inst. Hautes \'Etudes Sci. Publ. Math.}, (47):269--331 (1978),
  1977.

\end{thebibliography}
%\bibliographystyle{abbrv}

%\pagebreak \

\vfill

\begin{center}
\noindent
\begin{minipage}{\linewidth}
\small \noindent \textsc
{Manuel Amann} \\
\textsc{Department of Mathematics}\\
\textsc{University of Toronto}\\
\textsc{Earth Sciences 2146}\\
\textsc{Toronto, Ontario}\\
\textsc{M5S 2E4} \\
\textsc{Canada}\\
[1ex]
\textsf{mamann@uni-muenster.de}\\
\textsf{http://individual.utoronto.ca/mamann/}
\end{minipage}
\end{center}

\end{document}